\documentclass[11pt]{amsart}
\usepackage{amsmath, amsthm, amssymb, amscd, amsfonts}
\usepackage[all]{xy}

\def\pf{\begin{proof}}
\def\epf{\end{proof}}

\newcommand{\End}{\operatorname{End}}
\newcommand{\alga}{\Aut_{\rm Alg} A}
\newcommand{\algo}{(\alga)_0}

\newcommand{\Ss}{{\mathcal S}}
\newcommand{\hg}{{\mathcal H}}
\newcommand{\tg}{{\mathcal T}}
\newcommand{\ha}{{\mathbb H}}

\newcommand\ev{\operatorname{ev}}
\newcommand\coev{\operatorname{coev}}
\newcommand{\Aut}{\mbox{\rm Aut\,}}

\newcommand\id{\operatorname{id}}

\newcommand{\Ima}{\mbox{\rm Im\,}}
\newcommand{\spec}{\mbox{\rm Spec\,}}
\newcommand{\Ker}{\mbox{\rm Ker\,}}

\newcommand{\ydhh}{{}_H^{H}\mathcal{YD}}
\newcommand{\Autdiag}{\mbox{\rm Autdiagr\,}}
\newcommand{\ku}{{\Bbbk}}
\newcommand{\toba}{{\mathfrak B}}
\newcommand{\Z}{{\mathbb Z}}
\newcommand{\N}{{\mathbb N}}

\numberwithin{equation}{section}
\theoremstyle{plain}
\newtheorem{theo}{Theorem}[section]
\newtheorem{lemma}[theo]{Lemma}
\newtheorem{corollary}[theo]{Corollary}
\newtheorem{proposition}[theo]{Proposition}

\theoremstyle{definition}
\newtheorem{defi}[theo]{Definition}
\newtheorem{obs}[theo]{Remark}

\newtheorem{problem}{Problem}

\begin{document}

\title{On the automorphisms of $U_q^+({\mathfrak g})$}
\author{Nicol\'as Andruskiewitsch and Fran\c cois Dumas}
\address{Facultad de Matem\'atica, Astronom\'\i a y F\'\i sica
\newline \indent
Universidad Nacional de C\'ordoba
\newline
\indent CIEM -- CONICET
\newline
\indent (5000) Ciudad Universitaria, C\'ordoba, Argentina}
\email{andrus@mate.uncor.edu, \quad \emph{URL:}\/
http://www.mate.uncor.edu/andrus}
\address{Universit\'e Blaise Pascal
\newline \indent
Laboratoire de Math\'ematiques
(UMR 6620 du CNRS)
\newline \indent
F-63177, Aubi\`ere (France)}
\email{Francois.Dumas@math.univ-bpclermont.fr}
\date{\today}
\thanks{ 
\newline \indent This work was partially supported by Agencia C\'ordoba
Ciencia, ANPCyT-Foncyt, CONICET, ECOS, International Project
CNRS-CONICET ``M\'etodos hom\'ologicos en representaciones y
\'algebras de Hopf", PICS-CNRS 1514 and Secyt(UNC)}

\begin{abstract}
Let ${\mathfrak g}$ be a simple complex finite dimensional Lie
algebra and let $U_q^+({\mathfrak g})$ be the positive part of the
quantum enveloping algebra of $\mathfrak g$. If $\mathfrak g$ is
of type $A_2$, the group of algebra automorphisms of
$U_q^+({\mathfrak g})$ is a semidirect product $(\ku^{\times})^2
\rtimes \Autdiag ({\mathfrak g})$; any algebra automorphism is an
automorphism of braided Hopf algebra, and preserves the standard
grading \cite{AD1}. This intriguing smallness of the group of
algebra automorphisms raises questions about the extent of these
phenomena. We discuss some of them in the present paper. We
introduce the notion of ``algebra with few automorphisms" and
establish some consequences. We prove some exploratory results
concerning the group of algebra automorphisms for the type $B_2$.
We study Hopf algebra automorphisms of Nichols algebras and their
bosonizations and compute in particular the group of Hopf algebra
automorphisms of $U_q^+({\mathfrak g})$.
\end{abstract}

\maketitle
\section*{Introduction}
Let ${\mathfrak g}$ be a finite dimensional  complex simple Lie
algebra and ${\mathfrak g}={\mathfrak n}^- \oplus{\mathfrak
h}\oplus {\mathfrak n}^+$ be a triangular decomposition of
${\mathfrak g}$ related to a Cartan subalgebra ${\mathfrak h}$.
The structure of the group $\Aut_{\rm Alg} U_q({\mathfrak g})$ of
algebra automorphisms of the quantum enveloping algebra
$U_q({\mathfrak g})$ seems to be known only in the elementary case
where ${\mathfrak g}$ is of type $A_1$ (see \cite{AC} or
\cite{AD1}). The automorphism group $\Aut_{\rm Alg} \check
U_q({\mathfrak b}^+)$ of the augmented form of the quantum
enveloping algebra of the Borel subalgebra $\mathfrak
b^+={\mathfrak h}\oplus {\mathfrak n}^+$ is described for any
$\mathfrak g$ in \cite{Fl}. The groups of Hopf algebra
automorphisms of $U_q({\mathfrak g})$ and $\check U_q({\mathfrak
b}^+)$ are determined in \cite{CM} and \cite{Fl} respectively. We
are concerned in this paper with the group of automorphisms of the
quantum enveloping algebra $U_q^+({\mathfrak g})=U_q({\mathfrak
n}^+)$ of the nilpotent part.

\medbreak  Let $(V,c)$ be a braided vector space, see \ref{bvs}
below, and let  $\toba(V)$ be the corresponding Nichols algebra.
We prove that the group $\Aut_{\rm Hopf}\,\toba(V)$ of braided
Hopf algebra automorphisms of $\toba(V)$ coincides with the group
$GL(V,c)$ of automorphisms of braided vector space of $(V,c)$.
Thus we have
$$
\Aut_{\rm Hopf}\,\toba(V) \subset
\Aut_{\rm Gr Alg}\,\toba(V) \subset \Aut_{\rm Alg}\,\toba(V),
$$
where $\Aut_{\rm Gr Alg}$ means the group of algebra automorphisms
homogeneous of degree 0, in this case with respect to the standard
grading of $\toba(V)$, and $\Aut_{\rm Alg}$ is the group of all
algebra automorphisms.

\medbreak The class of Nichols algebras includes symmetric
algebras, free algebras, Grassmann algebras; thus there is no hope
to have a round description of the group $\Aut_{\rm
Alg}\,\toba(V)$. In particular the computation of this group is a
well-known classical open problem in the case of the symmetric
algebras.

\medbreak Let $\ku$ be the ground field and let $(V,c)$ be a
braided vector space of diagonal type, $n = \dim V$. Under some
technical assumptions, the group $GL(V,c)$ reduces in this case to
the semidirect product of the canonical action of the torus
$(\ku^{\times})^{n}$  by the subgroup $\Autdiag (c)$ of the
symmetric group $\mathbb S_n$ preserving the matrix of the
braiding.

 \medbreak A fundamental characterization due to Lusztig and Rosso
(in different but equivalent formulations), says that
$U_q^+({\mathfrak g}) = \toba(\mathfrak h)$ with a diagonal
braiding $c$ with matrix $(q^{d_ia_{ij}})$ build up from $q$ and
the Cartan matrix of $\mathfrak g$; say $n = \dim \mathfrak h$.
Furthermore, $\Autdiag c$ is in the case the group
$\Autdiag(\mathfrak g)$ of automorphisms of the Dynkin diagram of
$\mathfrak g$. Therefore, \bigbreak

\centerline{$\Aut_{\rm Hopf}\,U_q^+({\mathfrak
g})\simeq(\ku^{\times})^{n} \rtimes \Autdiag ({\mathfrak
g})$.}\bigbreak

The group $\Aut_{\rm Alg} U_q^+({\mathfrak g})$ was not determined
yet, to the best of our knowledge. However, if $\mathfrak g$ of
type $A_2$, then $\Aut_{\rm Alg} U_q^+({\mathfrak
g})\simeq(\ku^{\times})^{2} \rtimes \Autdiag ({\mathfrak g})$
\cite{AD1} (see also \cite{C2}). Thus, $\Aut_{\rm Hopf}\,\toba(V)
= \Aut_{\rm Gr Alg}\,\toba(V) = \Aut_{\rm Alg}\,\toba(V)$ in this
case. This intriguing result motivates several questions:

\begin{problem}\label{pb} Is it true that $\Aut_{\rm Alg} U_q^+({\mathfrak
g})\simeq(\ku^{\times})^{n} \rtimes \Autdiag ({\mathfrak g})$
?\end{problem}

We conjecture that the answer is positive for any $\mathfrak g$.
But even if the answer were negative, we would still ask: is it
true that $\Aut_{\rm Alg} U_q^+({\mathfrak g})$ is an algebraic
group?

\begin{problem}\label{pb2} Determine the braided vector spaces $(V, c)$
such that $$\Aut_{\rm Hopf}\,\toba(V) = \Aut_{\rm Gr
Alg}\,\toba(V) = \Aut_{\rm Alg}\,\toba(V).$$ \end{problem}

Graded algebras $A$ with the property $\Aut_{\rm Gr Alg}\,A =
\Aut_{\rm Alg}\, A$ do not seem to abound. Thus, we dare to pose:

\begin{problem}\label{pb3} Classify graded algebras $A$ with the
property $\Aut_{\rm Gr Alg}\,A = \Aut_{\rm Alg}\, A$.
\end{problem}

In this paper we contribute mainly to the first Problem. Let us
review the contents of the article. The first section is about
Hopf algebra automorphisms of Nichols algebras and their
bosonizations. We obtain from some general considerations the
computation of the group $\Aut_{\rm Hopf}\,U_q(\mathfrak b^+)
\simeq(\ku^{\times})^{n} \rtimes \Autdiag({\mathfrak g})$,
recovering a result from \cite{Fl}.

\bigbreak In section two, we briefly recall the case of type
$A_2$. We introduce the notion of ``algebras with few
automorphisms"; we classify gradings of these algebras and apply
to $U_q^+({\mathfrak{sl}_3})$.

\bigbreak Section three is an exploration of the case where
$\mathfrak g$ is of type $B_2$. Since there are no nontrivial
diagram automorphism in this case, the question is then to
determine if $\Aut_{\rm Alg}\,U_q^+({\mathfrak g})$ is isomorphic
or not to $(\ku^{\times})^2$.

\bigbreak A basic idea to approach the group $\Aut_{\rm
Alg}\,U_q^+({\mathfrak g})$ is to study its actions on natural
sets. In the paper \cite{AD1}, the study of the actions on the
sets of central and normal elements was crucial. This method fails
here because the center of $U_q^+({\mathfrak g})$ is a polynomial
algebra $\ku[z,z']$ in two variables (which are homogeneous
elements of degree 3 and 4 respectively for the canonical grading)
and any normal element of $U_q^+({\mathfrak g})$ is automatically
central.

\bigbreak The next natural sets where our group acts are the
various Spectra; the investigation of these actions is the matter
of this section. We begin by the ideals $(z)$ and $(z')$; these
are completely prime of height one and the factor domain
$U_q^+({\mathfrak g})/(z)$ is isomorphic to the quantum Heisenberg
algebra $U_q^+({\mathfrak sl}_3)$. Using the results of section
two, this allows to separate up to isomorphism the factor domains
$U_q^+({\mathfrak g})/(z)$ and $U_q^+({\mathfrak g})/(z')$. We
then show, first, that automorphisms of $U_q^+({\mathfrak g})$
cannot exchange the ideals $(z)$ and $(z')$; and second, that the
subgroup of $\Aut_{\rm Alg}\,U_q^+({\mathfrak g})$ of
automorphisms preserving the ideal $(z)$ reduces to the torus
$(\ku^{\times})^2$. To progress further in this direction, we need
better knowledge on the prime ideals of height one. Although the
stratification of the prime spectrum is known \cite{G}, the full
classification of the prime ideals is still open. We discuss the
stratification for type $B_2$ in subsection \ref{hspec}.

\bigbreak

\subsection*{Convention} We denote by $\N$ the set of non-negative
integers $\{0, 1, 2, 3, \dots \}$

\subsection*{Acknowledgements}
The origin of this paper lies in conversations with Jacques Alev,
during visits of the first author to the University of Rheims. His
enthusiasm about the automorphism problem was decisive to convince
us to consider this question. We also thank G\'erard Cauchon and
St\'ephane Launois for many enlightening discussions about the
third part of this paper.

This joint work was partially realized during visits of the first
author to the University of Clermont-Ferrand in March 2002 and
June 2003, and a visit of the second author to the University of
Cordoba in November 2003, in the framework of the Project ECOS
conducted by J.-L. Loday and M. Ronco, of the Project CONICET-CNRS
"M\'etodos hom\'ologicos en representaciones y \'algebras de Hopf"
and of the Project PICS-CNRS 1514 conducted by C. Cibils.

\bigbreak
\section{Braided Hopf algebra automorphisms}\label{nich}

In subsections 1.3 and 1.2 the field $\ku$ is arbitrary; in 1.3,
$\ku$ has characteristic 0 and contains an element $q$ not
algebraic over $\mathbb Q$.

\subsection{Braided vector spaces}

\subsubsection{Braided vector spaces}\label{bvs}
A braided vector space is a pair $(V, c)$ where $V$ is a vector
space $V$ over $\ku$ and $c: V\otimes V\to V\otimes V$ is a linear
isomorphism that is a solution of the braid equation $(c\otimes
\id) (\id\otimes c) (c\otimes \id) = (\id\otimes c) (c\otimes \id)
(\id\otimes c)$.

\bigbreak A  braided vector space $(V, c)$ is \emph{rigid} if $V$
is finite dimensional and the map $c^{\flat}:V^{*}\otimes V\to
V\otimes V^{*}$ is invertible, where
$$c^{\flat}=(\ev_V\otimes\id_{V\otimes V^{*}}) (\id_{V^{*}}\otimes
c\otimes\id_{V^{*}}) (\id_{V^{*}\otimes V}\otimes \coev_V).$$
Here $\ev_V: V^{*}\otimes V\to \ku$ is the usual evaluation map,
and $\coev_V: \ku\to V\otimes V^{*}$ is the coevaluation (the transpose of the trace).

\bigbreak
\subsubsection{Automorphisms of a braided vector space of diagonal type}\label{autbvs}
Let $(V, c)$ be a braided vector space. The braiding $c: V\otimes
V\to V\otimes V$ is said to be {\it diagonal} if there exists a
basis $x_{1}, \dots, x_{n}$ of $V$ and a matrix $(q_{ij})_{1\leq
i,j\leq n}$ with entries in $\ku^{\times}$ such that
$c(x_{i}\otimes x_{j}) = q_{ij} x_{j}\otimes x_{i}$ for any $1\leq
i,j\leq n$. In particular, $V$ is rigid and has finite dimension
$n\geq 1$.

\begin{obs} If the braiding is diagonal, then the matrix $(q_{ij})_{1\leq
i,j\leq n}$ does not depend on the basis $x_{1}, \dots, x_{n}$, up
to permutation of the index set $\{1, \dots, n\}$, see \cite[Lemma
1.2]{as-crelle}.
\end{obs}

\bigbreak Let $(V, c)$ be a braided vector space. A linear
automorphism $g\in GL(V)$ is said to be a \emph{braided vector
space automorphism} of $(V,c)$ if $g\otimes g$ commutes with $c$.
We denote by $GL(V,c)$ the corresponding subgroup of $GL(V)$.

\bigbreak Suppose that $c$ is of diagonal type for some basis
$x_{1}, \dots, x_{n}$ of $V$ and some matrix $(q_{ij})_{1\leq
i,j\leq n}$. We consider the following subgroup of the symmetric
group $\mathbb S_{n}$: \bigbreak

\centerline{ $\Autdiag (c) : =\{\sigma \in \mathbb S_{n} : q_{ij}
= q_{\sigma(i), \sigma (j)}, \ {1\le i, j \le n}\}$.}\bigbreak

Any $\sigma$ in $\Autdiag (c)$ induces naturally an automorphism
$g_{\sigma}\in GL(V,c)$ by $g_{\sigma}(x_j)=x_{\sigma(j)}$
for $1\leq j\leq n$. Any $g\in GL(V,c)$ of this type is called a
\emph{diagram automorphism} of $(V,c)$. Moreover it is clear that
the torus $(\ku^{\times})^n$ acts on $V$ by braided vector space
automorphisms. The following lemma gives necessary conditions for
the group $GL(V,c)$ to be generated by these two particular
subgroups.
\medbreak

\begin{lemma}\label{autbvsdiag}
Let $(V, c)$ be a braided vector space of diagonal type, with
respect to a basis $x_{1}, \dots, x_{n}$ of $V$ and a matrix
$(q_{ij})_{1\leq i,j\leq n}$ with entries in $\ku^{\times}$.
Assume that at least one of the following
conditions is satisfied:

\begin{itemize}
\item[(i)] For any $i \neq j$, there exists $h$ such that $q_{ih} \neq q_{jh}$.

\item[(ii)] For any $i \neq j$, there exists $h$ such that $q_{hi} \neq q_{hj}$.

\item[(iii)] For any $i \neq j$, the matrix  $\left(\begin{smallmatrix}
q_{ii} & q_{ij}\\ q_{ji} & q_{jj}\end{smallmatrix}\right)$ is not of the form
 $\left(\begin{smallmatrix} q & q\\ q & q\end{smallmatrix}\right)$.
\end{itemize}

Then we have

\centerline{$GL(V,c) \simeq (\ku^{\times})^{n} \rtimes \Autdiag (c)$.}
\end{lemma}

\pf Let $g\in GL(V)$; denote $g(x_i)=\sum_s\lambda_{s, i} x_s$,
$1\le i \le n$. Then

\begin{align*}
(g \otimes g)(c(x_i \otimes x_j))
&= (g \otimes g)( q_{ij} (x_j \otimes x_i))
= \sum_{1\le r,s \le n}q_{ij}\lambda_{r, j} \lambda_{s, i} \,
x_r \otimes x_s,
\\ c (g \otimes g)  (x_i \otimes x_j)  &= c\, \Bigl(
\sum_{1\le r,s \le n} \lambda_{r, j} \lambda_{s, i} \, x_s \otimes
x_r  \Bigr) = \sum_{1\le i,j \le n} q_{sr} \lambda_{r, j}
\lambda_{s, i} \, x_r \otimes x_s.
\end{align*}

Therefore $g\in GL(V,c)$ is and only if
\begin{equation}\label{tres}
q_{ij}\lambda_{r, j} \lambda_{s, i} = q_{sr}
\lambda_{r, j} \lambda_{s, i}, \quad \text{for all}\ 1 \le i,j,r,s
\le n.
\end{equation}

Suppose $g\in GL(V,c)$. Since $g$ is invertible, there exists
$\sigma \in \mathbb S_{n}$ such that $\lambda_{\sigma(h), h} \neq
0$ for any $1\le h \le n$. Then we deduce from \eqref{tres} that $
q_{ij} = q_{\sigma(i), \sigma (j)}$ for all $1\le i, j \le n$.

\bigbreak
Assume first that (i) is satisfied, and choose
$i, s \in \{1, \dots, n\}$ such that $\lambda_{s, i} \neq 0$.
Apply \eqref{tres} with any $j$ and $r=\sigma(j)$.
We obtain $q_{s, \sigma (j)} = q_{ij} = q_{\sigma(i), \sigma (j)}$,
for any $1 \le j \le n$; then $s = \sigma(i)$.
This implies that $g(x_i)=\lambda_{\sigma(i), i} x_{\sigma(i)}$ for any
$1\le i \le n$, and proves the result.

\bigbreak
Assume now that (ii) is satisfied, and choose
$j,r \in \{1, \dots, n\}$ such that $\lambda_{r,j} \neq 0$.
Apply \eqref{tres} with any $i$ and $s=\sigma(i)$.
We obtain $q_{\sigma (i), r} = q_{ij} = q_{\sigma(i), \sigma (j)}$,
for any $1 \le i \le n$; then $r = \sigma(j)$ and we conclude as in the previous case.

\bigbreak Finally assume that we have (iii) and take $i = j$ in \eqref{tres}.
If $\lambda_{s, i} \neq 0$ then $q_{ii} = q_{\sigma (i), s}
= q_{s,\sigma(i)} = q_{\sigma(i), \sigma (i)}$; therefore $s = \sigma(i)$
and we conclude as above.\epf

\bigbreak
\subsubsection{Braided Hopf algebras}
A non-categorical version of the concept of braided Hopf algebra was studied in \cite{Tk}.
A braided Hopf algebra is a collection $(R, m, \Delta, c)$ such that
\begin{itemize}
\item $(R, c)$ is a braided vector space,
\item $(R, m)$ is an associative algebra with unit $1$,
\item $(R, \Delta)$ is a coassociative coalgebra with counit $\varepsilon$,
\item $m, \Delta, 1, \varepsilon$ commute with $c$ in the sense of \cite{Tk},
\item $\Delta \circ m = (m\otimes m) (\id \otimes c\otimes \id)(\Delta \otimes \Delta)$,
\item the identity has an inverse for the convolution product in
$\End R$ (this inverse is called the \emph{antipode} and denoted by $\Ss$).
\end{itemize}
Here, recall that the convolution product of $f, g\in \End R$ is
given by $f * g = m(f\otimes g)\Delta$.

\bigbreak A homomorphism of braided Hopf algebras is a linear map
preserving $m, \Delta, c$.

\medbreak
\begin{lemma}\label{schau}
Let $R$ be a braided Hopf algebra and let $T: R \to R$
be a linear isomorphism that is an algebra and coalgebra map.
Then $T$ is a morphism of braided Hopf algebras.
\end{lemma}

\pf Let us define $T.f := TfT^{-1}$, $f\in \End R$. Then
$$
T.(f * g) = T \left(m(f\otimes g)\Delta\right)T^{-1} =
m(T\otimes T)(f\otimes g)(T^{-1}\otimes T^{-1})\Delta
= (T.f) * (T.g);$$
hence $T.\Ss = \Ss$, or $T\Ss = \Ss T$. But it was shown in \cite{Sb}
that the braiding of a braided Hopf algebra can be expressed in
terms of the product, coproduct and antipode. Thus $T$ preserves also the braiding $c$.
\epf

The group of braided Hopf algebra automorphisms of $R$ is denoted by
$\Aut_{\rm Hopf}\,R$.

\bigbreak
\subsubsection{Yetter-Drinfeld modules}\label{ydmod}
Yetter-Drinfeld modules give rise to braided vector
spaces and play a fundamental r\^ole in problems related to the
classification of Hopf algebras.

\bigbreak Let us recall that a Yetter-Drinfeld module $V$ over a
Hopf algebra $H$ with bijective antipode $\Ss$ is both a left
$H$-module and left $H$-comodule, such that the action $H\otimes
V\dot\rightarrow V$ and the coaction $\delta:V\rightarrow H\otimes V$
satisfy the compatibility condition:
$\delta(h.v)=h_{(1)}v_{(-1)}\Ss h_{(3)}\otimes h_{(2)}.v_{(0)}$ for
all $h\in H,v\in V$.
We denote by $\ydhh$ the category of Yetter-Drinfeld
modules over $H$, where morphisms respect both the action and
the coaction of $H$.

\bigbreak The usual tensor product defines a structure of monoidal
category on $\ydhh$. It is braided, with braiding $c_{V, W}:
V\otimes W \to W\otimes V$ defined by $c(v \otimes w) = v_{(-1)}.w
\otimes v_{(0)}$, for $V, W \in \ydhh$, $v\in V$, $w\in W$. Then
$(V,c_{V, V})$ is a braided vector space, for any $V\in \ydhh$. It
is known that any rigid  braided vector space can be realized as a
Yetter-Drinfeld module over a (non-unique) Hopf algebra,
essentially by the FRT-construction; see \cite{Tk} for references
and details.

\bigbreak As in any braided monoidal category, there is the notion
of Hopf algebras in $\ydhh$. Hopf algebras in $\ydhh$ are braided
Hopf algebras by forgetting the action and the coaction.
Conversely, let $R$ be any braided Hopf algebra whose underlying
braided vector space is rigid. Then there exists a (non-unique)
Hopf algebra $H$ such that $R$ can be realized as a Hopf algebra
in $\ydhh$ \cite{Tk}.

\bigbreak
\subsubsection{Bosonizations of a Hopf algebra in $\ydhh$}
\label{radford}

We recall the bosonization procedure, or Radford biproduct, found
by Radford and explained in terms of braided categories by Majid.

\bigbreak
Let $H$ be a Hopf algebra with bijective antipode. Let $R$ be a Hopf algebra in $\ydhh$.
The bosonization of $R$ by $H$ is the (usual) Hopf algebra $A=R\#H$, with underlying
vector space $R\otimes H$, whose multiplication and comultiplication are
given by:\medbreak

\centerline{
$(r\#h)(s\#f)=r(h_{(1)}.s)\#h_{(2)}f$ \quad  and \quad
$\Delta(r\#h)=r^{(1)}\#(r^{(2)})_{(-1)}h_{(1)}\otimes (r^{(2)})_{(0)}\#h_{(2)}$.}\medbreak

The maps $\pi:A\rightarrow H, r\#h\mapsto \epsilon(r)h$ and
$\iota : H\rightarrow A, h\mapsto 1\#h$ are Hopf algebra homomorphisms
and $R=\{a\in A: (\id\otimes\pi)\Delta(a)=a\otimes 1\}$.

\bigbreak Conversely, let $A$, $H$ be Hopf algebras with bijective
antipode and let $\pi:A\rightarrow H$ and $\iota : H\rightarrow A$
be Hopf algebra homomorphisms such that $\pi \iota = \id_H$. Then
$R=\{a\in A: (\id\otimes\pi)\Delta(a)=a\otimes 1\}$ is a Hopf
algebra in $\ydhh$ and the multiplication induces an isomorphism
of Hopf algebras $R\# H \simeq A$.

\bigbreak
\subsubsection{Nichols algebras} We recall the
definition of Nichols algebra, see \cite{AS} for details and
references.

Let $V\in \ydhh$. A graded Hopf algebra $R= \bigoplus_{n\ge 0}
R(n)$ in $\ydhh$ is called a {\it Nichols algebra} of $V$ if
$\ku\simeq R(0)$ and $V\simeq R(1)$ in $\ydhh$, and if:
\begin{flalign} \label{nichols1}
&R(1) = \mathcal P(R), \text{ the space of primitive elements of }
R,&
\\ \label{nichols2}
&R \text{ is generated as an algebra by } R(1).&
\end{flalign}

The Nichols algebra of $V$ exists and is unique up to isomorphism;
it is denoted by $\toba(V)= \bigoplus_{n\ge 0} \toba^n(V)$. The
associated braided Hopf algebra (forgetting the action and the
coaction) depends only on the braided vector space $(V, c)$. We
shall identify  $V$ with the subspace of homogeneous elements of
degree one in $\toba(V)$.

\medbreak Let us recall the following explicit construction of
$\toba(V)$. For any integer $m\geq 2$, we denote by $\mathbb B_m$
the $m$-braid group. A presentation of $\mathbb B_m$ is given by
generators
 $\sigma_1,\ldots,\sigma_{m-1}$ and relations
 $\sigma_i\sigma_j=\sigma_j\sigma_i$ if $|i-j|\geq 2$
an $\sigma_i\sigma_{i+1}\sigma_i=\sigma_{i+1}\sigma_i\sigma_{i+1}$
for any $1\leq i\leq m-2$.
There is a natural projection $\pi: \mathbb B_{m} \to \mathbb S_{m}$
sending $\sigma_{i}$ to the transposition $\tau_{i} :=(i, i+1)$ for all $i$.
This projection $\pi$ admits a set-theoretical section
$s:\mathbb S_{m} \to \mathbb B_{m}$ determined by
\begin{align*}
s(\tau_{i}) = \sigma_{i}, \quad 1\le i \le n-1,\qquad
s(\tau\omega) = s(\tau)s(\omega), \quad \text{if \ } \ell(\tau\omega) = \ell(\tau) + \ell(\omega).
\end{align*}
Here $\ell$ denotes the length of an element of $\mathbb S_{m}$
with respect to the set of generators $\tau_{1}, \dots, \tau_{m-1}$. The map $s$ is called the Matsumoto section.
In other words, if $\omega = \tau_{i_1} \dots \tau_{i_j}$
is a reduced expression of $\omega \in \mathbb S_{m}$,
then $s(\omega) = \sigma_{i_1} \dots \sigma_{i_j}$.
Using the section $s$, the following distinguished elements of the group algebra
$\ku\mathbb B_{m}$ are defined:
$$\mathfrak S_{m} := \sum_{\sigma \in \mathbb S_{m}} s(\sigma).
$$
By convention, we still denote by
$\mathfrak S_{m}$ the images of these elements in
$\End (T^{m}(V))= \End(V^{\otimes m})$ via the representation
$\rho_m:{\mathbb B}_m\rightarrow \Aut(V^{\otimes m})$ defined by
$\rho_m(\sigma_i)=\id\otimes\ldots\otimes\id\otimes c\otimes\id\otimes\ldots\otimes\id$,
with $c$ acting on the tensor product of the copies of $V$ indexed by $i$ and $i+1$.

\bigbreak Let $\widetilde c$ be the canonical extension of $c$
into a braiding of $T(V)$. Let $T(V)\,\underline{\otimes}\,
T(V)$ be the algebra whose underlying vector space is $T(V)\otimes
T(V)$ with the product ``twisted" by  $\widetilde c$. There is a
unique algebra map $\Delta: T(V)\to T(V)\,\underline{\otimes}\,
T(V)$ such that $\Delta(v) v\otimes 1 + 1 \otimes v$, $v\in V$.
Then $T(V)$ is a braided Hopf algebra and $\toba (V) = T(V)/J$
where $J=\bigoplus_{m \ge 0}\Ker\mathfrak S_{m}$ (see for instance
\cite{AS}, \cite{N} or \cite{Ro}). \medbreak

\subsection{Automorphisms of Nichols algebras and their bosonizations}

\subsubsection{Hopf algebra automorphisms of a Nichols algebra}
We can now compute the group of Hopf algebra automorphisms of a Nichols algebra, \emph{cf.} the notation introduced in \ref{autbvs}.

\medbreak
\begin{theo}\label{authopfnich}
There is a group isomorphism
$\toba: GL(V,c) \to \Aut_{\rm Hopf}\,\toba(V)$.\end{theo}

\pf For any $g\in GL(V)$, we denote by $\widetilde g$ the
canonical extension of $g$ into an algebra automorphism of $T(V)$.
If $g$ commutes with $c$, then $\widetilde g$ commutes with
$\sum_{m \ge 2} \mathfrak S_{m}$, and so induces an algebra
automorphism $\toba(g)$ of $\toba(V)=T(V)/J$. In order to prove that $\toba(g)$ is also a coalgebra automorphism of $\toba (V)$,
we claim that $\widetilde g$ is a coalgebra map: $((\widetilde g\otimes\widetilde g)\circ\Delta)(v)=
(\Delta\circ\widetilde  g)(v)$, for any $v\in T(V)$.

It is clear that this assertion is true when $v\in V$. So it is
enough to prove that $\widetilde g\otimes\widetilde g$ is a
morphism of the algebra $T(V)\,\underline{\otimes}\,  T(V)$. For that, let us consider
$u,v,x,y\in T(V)$. Let us set $\widetilde c(y\otimes u)=u'\otimes y'$ in a symbolic way. By definition of the twisted product
in $T(V)$, with have $(x\otimes y)(u\otimes v)=xu'\otimes y'v$,
then: \medbreak

\centerline{ $(\widetilde g\otimes\widetilde g)((x\otimes
y)(u\otimes v))= \widetilde g(x)\widetilde g(u')\otimes \widetilde
g(y')\widetilde g(v)$.}\medbreak

Moreover, it is easy to check that the assumption $g\in GL(V,c)$
implies $\widetilde g\in GL(T(V),\widetilde c)$, so $(\widetilde
g\otimes\widetilde g)(u'\otimes y')= \widetilde c(\widetilde
g(y)\otimes\widetilde g(u))$, and then: \medbreak

\centerline{ $(\widetilde g\otimes\widetilde g)(x\otimes y)
(\widetilde g\otimes\widetilde g)(u\otimes v) = (\widetilde
g(x)\otimes\widetilde g(y)) (\widetilde g(u)\otimes\widetilde
g(v)) = \widetilde g(x)\widetilde g(u')\otimes\widetilde
g(y')\widetilde g(v)$.}\medbreak

Hence $(\widetilde g\otimes\widetilde g)((x\otimes
y)(u\otimes v))= (\widetilde g\otimes\widetilde g)(x\otimes y)
(\widetilde g\otimes\widetilde g)(u\otimes v)$, as claimed.
By Lemma \ref{schau}, $\toba(g)$ preserves $c$.
Thus, we have a well-defined map $\toba: GL(V,c) \to
\Aut_{\rm Hopf}\,\toba(V)$, which is injective by \eqref{nichols2}.

\bigbreak
Conversely, let $u: \toba(V) \to \toba(V)$ be an
automorphism of braided Hopf algebras.
Then $u(V) = V$ since $V = P(\toba(V))$ and the theorem follows. \epf

\bigbreak
\subsubsection{Automorphisms of bosonizations}

Our goal is to compute the Hopf algebra automorphisms of $A=R\#H$
when $R$ is a Nichols algebra, under suitable hypothesis. The
exposition is inspired by \cite[Section 6]{as-adv}. We begin by a
description of a natural class of such automorphisms, for general
$R$.

\medbreak
\begin{lemma}\label{autboslema} Let $H$ be a Hopf algebra and let $R$ be a Hopf algebra in $\ydhh$.  Let $G: R \to R$ and $T: H \to H$ be linear maps. Then $G\# T:= G\otimes T : R\# H \to R\# H$ is a Hopf algebra map if and only if the following conditions hold:

\begin{flalign}
& T \text { is a Hopf algebra automorphism of } H, &
\\& G \text { is a Hopf algebra automorphism of } R,&
\\ \label{compatible3}
& G(h.s) = T(h). G(s), \qquad s\in R, \, h\in H,&
\\\label{compatible4}
& \delta \circ G = (T \otimes G)\circ \delta.&
\end{flalign}
\end{lemma}

\pf Left to the reader. \epf

A pair $(G,T)$ as in the Lemma shall be called \emph{compatible}.

\medbreak
\begin{lemma}\label{autboslema2}
Let $H$ be a Hopf algebra and $V$ a Yetter-Drinfeld module over
$H$.

\bigbreak\begin{itemize}
\item[(i)] Assume that $H$ is cosemisimple, and that
the following hypothesis holds:\bigbreak {\begin{itemize}
\item[(H)]
the types of the isotypic components of $V\# H$ under the
adjoint
action of $H$ do not appear in the adjoint action of $H$ on
itself.
\end{itemize}}\bigbreak
\noindent Then any Hopf algebra automorphism of $\toba(V) \# H$ is
of the form $G\# T$, with $(G,T)$ compatible.

\bigbreak
\item[(ii)] If in addition $H$ is commutative, then {\rm (H)} is equivalent to:\bigbreak
{\begin{itemize}
\item[(H$^\prime$)]
the trivial representation does not appear as a subrepresentation of $V$.
\end{itemize}}
\end{itemize}
\end{lemma}

\pf (i). Let $\Phi$ be a Hopf algebra automorphism of $A =
\toba(V) \# H$. It is known that the coradical filtration of $A$
is $A_m = \bigoplus_{0\le n \le m} \toba^n(V) \# H$
\cite[1.7]{AS}. Since $\Phi$ is a  coalgebra map, it preserves the
coradical filtration. In particular, $\Phi(H) = H$ and $\Phi(H
\oplus V \# H) = H \oplus V \# H$. Let $T: H\to H$ be the
restriction of $\Phi$; it is an automorphism of Hopf algebras.
Also, $\Phi: H \oplus V \# H \to H \oplus V \# H$ preserves the
adjoint action of $H$. By hypothesis (H), $\Phi(V \# H) = V \# H$.
Since $\Phi$ is an algebra map, this implies that $\Phi(\toba^n(V)
\# H) = \toba^n(V) \# H$, by \eqref{nichols2}. Let $\pi: A \to H$
be the projection with kernel $\bigoplus_{n \ge 1} \toba^n(V) \#
H$; clearly, $\Phi\pi = \pi\Phi$. Hence $\Phi(\toba(V)) =
\toba(V)$, since $\toba(V) = \{v \in A: (\id \otimes \pi)\Delta
(v) = v\otimes 1\}$. Let $G: \toba(V) \to \toba(V)$ be the
restriction of $\Phi$. Since $\Phi$ is an algebra map, $\Phi = G\#
T$.  By Lemma \ref{autboslema}, the pair $(G,T)$ is compatible.

\bigbreak (ii). If $H$ is commutative, the adjoint action of $H$
on itself is trivial, and the isotypic components of the adjoint
action of $H$ on $V\# H$ are of the form $U\# H$, where $U$ runs
in the set of isotypic components of the adjoint action of $H$ on
$V$. This shows that (H) is equivalent to (H$^\prime$) in this
case. \epf

The hypothesis (H) is needed, as the following example shows. Let
$A = \ku [x, g, g^{-1}]$ be the tensor product of the polynomial
algebra in $x$ and the Laurent polynomial algebra in $g$. This is
a Hopf algebra with $x$ primitive and $g$ group-like. The Hopf
algebra automorphism $T: A\to A$, $T(g) = g$, $T(x) = x + 1-g$,
does not preserve the Nichols algebra $\ku[x]$.

\bigbreak We now consider the following particular setting. We
assume that $H = \ku \Gamma$ is the group algebra of an abelian
group. We also assume the existence of a basis $x_1, \dots, x_n$
of $V$ such that, for some elements $g_1,\ldots,g_n\in\Gamma$,
$\chi_1,\ldots,\chi_n\in\widehat\Gamma$, the action and coaction
of $\Gamma$ are given by
$$h.x_j=\chi_j(h)x_j,\qquad \delta(x_j)=g(j)\otimes x_j, \qquad 1\leq j\leq n.$$

\bigbreak
\begin{theo}\label{autbos}
Suppose further that
\begin{flalign} \label{nichols3}
&\chi_i \neq \varepsilon, \qquad \qquad\qquad  1\leq i\leq n,&
\\ \label{nichols4}
&(g_i, \chi_i) \neq (g_j, \chi_j), \quad 1\leq i\neq j\leq n.&
\end{flalign}

Then there is a bijective correspondence between $\Aut_{\rm
Hopf}\,\toba(V) \# H$, and the set of pairs $(\phi, \psi)$, where
$\psi$ is a group automorphism of $\Gamma$ and $\phi: V \to V$ is
a linear isomorphism given by $\phi(x_i) =
\lambda_ix_{\sigma(i)}$, $1\leq i\leq n$, with $\lambda_i \in
\ku^{\times}$ and $\sigma \in \mathbb S_n$, such that
\begin{equation}\label{nichols5}
\psi(g_i) = g_{\sigma(i)}, \qquad \chi_i
= \chi_{\sigma(i)} \circ \psi, \qquad 1\leq i\leq n.\end{equation}
\end{theo}

\bigbreak \pf Condition \eqref{nichols3} guarantees that
hypothesis (H$^\prime$) holds. Thus any Hopf algebra automorphism
of $\toba(V) \# H$ is of the form $G\# T$, with $(G,T)$
compatible, by Lemma \ref{autboslema2}. But $T$ is determined by a
group automorphism $\psi$ of $\Gamma$, and $G$ is of the form
$\toba(\phi)$ for some $\phi \in GL(V,c)$ by Theorem
\ref{authopfnich}. Now conditions \eqref{compatible3} and
\eqref{compatible4} imply that
$$
\gamma. \phi(x_i) = (\chi_i \circ \psi^{-1})(\gamma) \phi(x_i), \qquad
\delta(\phi(x_i)) = \psi(g_i) \otimes \phi(x_i), \qquad \gamma \in \Gamma, \, 1 \le i \le n,$$
and thus
$$
\phi(x_i) \in \sum_{j\ : \ \psi(g_i)=g_j, \ \chi_i \circ \psi^{-1}
= \chi_j} \ku\, x_j, \qquad 1 \le i \le n.
$$
But condition \eqref{nichols4} implies that there is only one $j$
such that $\psi(g_i)=g_j$, $\chi_i \circ \psi^{-1} = \chi_j$; set
$\sigma(i) = j$. This defines $\sigma$, and clearly
\eqref{nichols5} holds.

Conversely, any pair $(\phi, \psi)$ as above gives raise to a
compatible pair $(G,T)$ with $G = \toba(\phi)$ and $T$ determined
by $\psi$. \epf

\subsection{Hopf algebra automorphisms of Nichols algebras of Drinfeld-Jimbo type}

\subsubsection{Definition and notations ({\it cf.}
\cite{BG},\cite{J},\cite{L})}\label{qea} We fix $q\in
\ku^{\times}$, $q$ not algebraic over $\mathbb Q$. Let ${\mathfrak
g}$ be a simple finite dimensional Lie algebra of rank $n$ over
$\ku$. Let ${\mathfrak g}={\mathfrak n}^-\oplus{\mathfrak
h}\oplus{\mathfrak n^+}$ be a triangular decomposition of
${\mathfrak g}$ related to a Cartan subalgebra $\mathfrak h$ of
$\mathfrak g$. Let $C=(a_{i,j})_{1\leq i,j\leq n}$ the associated
Cartan matrix and $(d_1,\ldots,d_n)$ the relatively primes
integers symmetrizing $C$. The quantum enveloping algebra of the
nilpotent positive part ${\mathfrak n}^+$ of $\mathfrak g$,
denoted by $U_q({\mathfrak n}^+)$ or $U_q^+({\mathfrak g})$, is
the algebra generated over $\ku$ by $n$ generators
$E_1,\ldots,E_n$ satisfying the quantum Serre relations: \medbreak

\centerline{ $\sum_{\nu=0}^{1-a_{i,j}} \left[\smallmatrix
1-a_{i,j}  \\ \nu  \\ \endsmallmatrix\right]_{q^{d_i}}
E_i^{1-a_{i,j}-\nu}E_j(-E_i)^{\nu}=0$ \quad for all $1\leq
i\not=j\leq n$.}\medbreak

The quantum enveloping algebra of the positive Borel algebra
$\mathfrak b^+=\mathfrak h\oplus\mathfrak n^+$, denoted by
$U_q({\mathfrak b}^+)$, is the algebra generated over $\ku$ by
$E_1,\ldots,E_n,K_1^{\pm 1},\ldots,K_n^{\pm 1}$ satisfying the
quantum Serre relations, the commutation between the $K_i$'s
and the $q$-commutation relations:\medbreak

\centerline{$K_iE_j=q^{d_ia_{i,j}}E_jK_i$ \quad for all $1\leq
i,j\leq n$.}
\medbreak

It is well-known that $U_q({\mathfrak b}^+)$ is a Hopf algebra for
the coproduct, counit and antipode defined by

\centerline{$\begin{matrix}
\Delta(K_i)=K_i\otimes K_i,\hfill & \varepsilon(K_i)=1, \hfill & S(K_i)=K_i^{-1},\hfill\\
\Delta(E_i)=E_i\otimes 1+ K_i\otimes E_i,\hfill & \varepsilon(E_i)=0, \hfill & S(E_i)=-K_iE_i.\hfill\\
\end{matrix}$}\medbreak

We denote by $H$ the Hopf subalgebra $H=\ku[K_1^{\pm
1},\ldots,K_n^{\pm 1}]\simeq \ku\Z^n$ of $U_q({\mathfrak
b}^+)$.

\bigbreak
\subsubsection{Braided Hopf algebra structure on $U_q^+({\mathfrak g})$}
From  corollary 33.1.5 of \cite{L} or theorem 15 of \cite{Ro}, we
have $U_q^+({\mathfrak g})=\toba(V)$ for $V={\ku}
E_1\oplus\ldots\oplus{\ku}E_n$ and $c$ the diagonal braiding of
$V$ defined from the Cartan matrix $C=(a_{i,j})_{1\leq i,j\leq n}$
and the integers $(d_1,\ldots,d_n)$ by

$$c(E_i\otimes E_j)=q^{d_ia_{i,j}}E_j\otimes
E_i.$$

Let $\iota: H \to U_q({\mathfrak b}^+)$ be the inclusion and let
$\pi: U_q({\mathfrak b}^+) \to H$ be the unique Hopf algebra map
such that $\pi(K_i)=K_i$, $\pi(E_i)=0$. Then $\pi \iota = \id_H$
and $U_q^+({\mathfrak g}) = \{a\in U_q({\mathfrak b}^+):
(\id\otimes\pi)\Delta(a)=a\otimes 1\}$. Hence $$U_q({\mathfrak
b}^+) \simeq U_q^+({\mathfrak g}) \# H,$$ \emph{cf.} Subsection
\eqref{radford}. Here the coaction is determined by $\delta(E_i) =
g_i \otimes E_i$, with $g_i = K_i$, $1\le i \le n$. Also, the
action is determined by $\gamma.E_i = \chi_i(\gamma) E_i$, for
$\gamma \in \Gamma = \Z^n$, where $\chi_i\in\widehat\Gamma$ is
defined by $\chi_i(K_j) = q^{d_ia_{i,j}}$,  $1\le i,j \le n$.

\bigbreak
\subsubsection{Automorphisms of $U_q^+({\mathfrak g})$ and $U_q^+({\mathfrak b})$}
With the notations of \ref{autbvs} for $q_{i,j}=q^{d_i}a_{i,j}$,
the subgroup $\Autdiag c$ of $\mathbb S_n$ is the group
$\Autdiag(\mathfrak g)$ of automorphisms of the Dynkin diagram of
$\mathfrak g$ (see for instance \cite{CM}), which acts by
automorphisms on $U_q^+({\mathfrak g})$ and $U_q({\mathfrak b}^+)$
by:\medbreak

\centerline{$\sigma\in\Autdiag(\mathfrak g): E_i\mapsto
E_{\sigma(i)}, \ K_i\mapsto K_{\sigma(i)}$ \quad for any $1\leq
i\leq n$.}\medbreak

The $n$-dimensional torus on $\ku$ also acts by automorphisms on
$U_q^+({\mathfrak g})$ and $U_q({\mathfrak b}^+)$, by:\medbreak

\centerline{$(\alpha_1,\ldots,\alpha_n)\in(\ku^{\times})^n: E_i\mapsto
\alpha_iE_i, \ K_i\mapsto K_i$ \quad for any $1\leq i\leq
n$.}\medbreak

Now we can prove the following theorem.

\medbreak
\begin{theo}\label{maintheo}
$\Aut_{\rm Hopf }\,U_q^+({\mathfrak g}) \simeq (\ku^{\times})^{n}
\rtimes \Autdiag ({\mathfrak g}) \simeq \Aut_{\rm Hopf
}\,U_q({\mathfrak b^+})$.\end{theo}

\pf The first isomorphism just follows from Lemma
\ref{autbvsdiag} and Theorem \ref{authopfnich}.

Let $\Gamma = \Z^n$, $g_i$ and $\chi_i$ as above. By  Theorem
\ref{autbos}, any $T \in \Aut_{\rm Hopf }\,U_q({\mathfrak b^+})$
is determined by a pair $(\phi, \psi)$, where $\psi$ is a group
automorphism of $\Gamma$ and $\phi: V \to V$ is a linear
isomorphism given by $\phi(x_i) = \lambda_ix_{\sigma(i)}$, $1\leq
i\leq n$, with $\lambda_i \in \ku^{\times}$ and $\sigma \in
\mathbb S_n$, such that \eqref{nichols5} holds. Then
$$
q^{d_ia_{i,j}} = \chi_i(K_j) =\chi_i\psi^{-1}\psi(K_j) =
\chi_{\sigma(i)}(K_{\sigma(j)}) = q^{d_{\sigma(i)}a_{\sigma(i),\sigma(j)}},
$$
hence $\sigma \in \Autdiag(\mathfrak g)$. Furthermore, $\psi$ is uniquely determined by $\sigma$.  This implies the second isomorphism.
\epf

\begin{obs}
The second isomorphism in the Theorem was proved previously in
\cite{Fl} as a corollary of the description of
the group of all algebra automorphisms
$\Aut_{\rm Alg}U_q({\mathfrak b^+})$ (in fact for the slightly
different augmented algebra $\check U_q({\mathfrak b}^+)$).
Effectively there exist automorphisms of the algebra $\check
U_q({\mathfrak b}^+)$ which are not Hopf algebra automorphisms; in
particular some combinatorial infinite subgroups of $n\times n$
matrices with coefficients in $\Z$, as well as the natural
action of the $2n$-dimensional torus on the $E_i$'s and $K_j$'s.
\end{obs}

\section{The case where $\mathfrak g$ is of the type $A_2$}

In this section the field $\ku$ has characteristic 0; in
subsection \ref{autheis}, $q\in \ku^{\times}$ is not a root of
one.

\subsection{Graded algebras with few automorphisms}\label{autfew}
Here we consider graded algebras with few automorphisms. To begin
with, we recall the well-known equivalence between gradings and
rational actions, see \cite[p. 150 ff.]{BG} and references
therein.

\bigbreak Let $\hg$ be an algebraic group. A representation $\rho:
\hg \to \Aut_{\ku} V$ of $\hg$ on a vector space $V$ is
\emph{rational} if  $V$ is union of finite dimensional rational
$\hg$-modules.

If $\hg = (\ku^{\times})^r$ is a torus, then there is a bijective
correspondence between
\begin{itemize}
\item rational actions of $\hg$ on $V$, and
\item gradings $V = \bigoplus_{m \in \Z^r} \, V_m$.
\end{itemize}
In this correspondence, $V_m$ is the isotypic component of type
$m$, where $\Z^r$ is identified with the group of rational
characters of $\hg$. Thus, $\hg$-submodules of $V$ are rational,
and they are exactly the graded subspaces of $V$.

\bigbreak Let $A$ be an associative algebra over $\ku$. A rational
action of $\hg$ on $A$ is one induced by a rational representation
$\rho: \hg \to \alga$ by algebra automorphisms. If $\hg =
(\ku^{\times})^r$ is a torus, then there is a bijective
correspondence between
\begin{itemize}
\item rational actions of $\hg$ on $A$, and
\item algebra gradings $A = \bigoplus_{m \in \Z^r} \, A_m$.
\end{itemize}

\begin{defi}\label{few}
An associative algebra $A$ has \emph{few automorphisms} if the
following conditions hold.
\begin{enumerate}
\item[(i)] There exists a finite dimensional
$\alga$-invariant subspace $V$ such that the restriction $\alga
\to GL(V)$ is injective; we identify $\alga$ with its image in
$GL(V)$;
\item[(ii)] $\alga$ is an algebraic subgroup of $GL(V)$;
\item[(iii)] the action of $\alga$ on $A$ is rational;
\item[(iv)] the connected component $\algo$ of the
identity of $\alga$ is isomorphic to a torus $(\ku^{\times})^r$.
\end{enumerate}
\end{defi}

\bigbreak
 Let $A$ be an algebra with few isomorphisms.
Then $A$ has a \emph{canonical} grading induced by the rational
action of $\algo \simeq (\ku^{\times})^r$:
$$A = \bigoplus_{m \in
\Z^r} \, A_{(m)}.$$

Let $\vert \, \vert: \Z^r \to \Z$ be the function $\vert m \vert =
\sum_{1\le j \le r} m_j$, if $m = (m_1, \dots, m_r)\in \Z^r$. The
$\Z$-grading induced by the canonical grading via $\vert \, \vert$
shall be called the \emph{standard grading} and denoted $A =
\bigoplus_{M \in \Z} \, A_{[M]}$. Thus
$$A_{[M]} = \bigoplus_{m \in
\Z^r : \vert m\vert = M} \, A_{(m)}, \qquad M\in \Z.$$

\bigbreak
\begin{lemma}
Let $A$ be an algebra with few automorphisms.\begin{itemize}

\item[(i)] Any algebra automorphism preserves the canonical grading, and
those in $\algo$ are homogeneous of degree 0.

\bigbreak
\item[(ii)] Assume that the canonical grading is nonnegative:
$A = \bigoplus_{m \in \N^r} \, A_{(m)}.$ Then any algebra
automorphism is homogeneous of degree 0 with respect to the
standard grading.
\end{itemize}
\end{lemma}

\pf (i). Let $\theta\in \alga$, let ${\rm inn}_{\theta}: \alga \to \alga$ be
the inner automorphism defined by $\theta$ and let $\widehat{{\rm
inn}_{\theta}}: \Z^r \to \Z^r$ be the induced group homomorphism. Then
$\theta\left(A_{(m)} \right) = A_{(\widehat{{\rm inn}_{\theta}}(m))}$ for any
$m\in \Z^r$, and this implies the claim.

\bigbreak (ii). Under this hypothesis, the matrix of
$\widehat{{\rm inn}_{\theta}}$ in the canonical basis has
non-negative entries. Thus $\theta\left(A_{[M]} \right) \subseteq
\sum_{s\ge 0}A_{[M] + s}$; but some power of $\theta$ belongs to
$\algo$, hence only $s = 0$ survives. \epf

Starting from the canonical grading, new gradings of $A$ can be
constructed by means of morphisms of groups $\Z^r \to \Z^t$; we
show now that no other grading arises in this case.

\begin{theo}\label{fewt}
Let $A$ be an algebra with few automorphisms and let
\begin{equation}\label{grading2}
A = \bigoplus_{n \in \Z^t} \, A_{n}
\end{equation}
be any algebra grading of $A$. Then there is a morphism of
groups $\varphi: \Z^r \to \Z^t$ such that
\begin{equation}\label{grading3}
A_n = \bigoplus_{m \in \Z^r: \varphi(m) = n} \, A_{(m)},
\end{equation}
for all $n\in \Z^t$.
\end{theo}

\pf Let $\tg$ be the torus $\Z^t$ and let $\rho$ be the
representation of $\tg$ induced by the grading \eqref{grading2}.
Let $V$ be the vector subspace as in Definition \ref{few}; since
$V$ is stable under $\alga$, it is also clearly stable under
$\tg$. Consider the commutative diagram
\begin{equation*}
\xymatrix{\ar[rd]_{\rho_{\vert V}} \tg  \ar[rr]^{\rho} & & \alga
\ar[dl]^{\rm res}
\\ &  GL(V).  & }
\end{equation*}
The map $\rho_{\vert V}$ is a homomorphism of algebraic groups;
then $\rho$ is homomorphism of algebraic groups, say by \cite[Ex.
3.10, p. 21]{H}. Since $\tg$ is connected, $\rho(\tg) \subseteq
\algo$. Thus, the transpose of $\rho$ induces a morphism of groups
$\varphi: \Z^r \to \Z^t$. But then $A_{(m)} \subset
A_{\varphi(m)}$. Since
$$A = \sum_{m \in \Z^r} \, A_{(m)} =
\sum_{n \in \Z^t}\Bigl(\, \sum_{m \in \Z^r: \varphi(m) = n} \,
A_{(m)}\Bigr) \subseteq \bigoplus_{n \in \Z^t} \, A_{n} = A,$$ we
get the equality \eqref{grading3}. \epf

\subsection{Algebra automorphisms of $U_q^+({\mathfrak{sl}_3})$}\label{autheis}

\subsubsection{Notations}\label{heis}

We suppose here that $\mathfrak g=\mathfrak{sl}_3$. Then we have
$n=2$, $C=\left(\begin{smallmatrix} 2 & -1\\ -1 & 2\\
\end{smallmatrix}\right)$, $d_1=d_2=1$, and $U_q^+(\mathfrak g)$
is the algebra generated over $\ku$ by $E_1$ and $E_2$ satisfying
the relations:
$E_1^2E_2-(q^2+q^{-2})E_1E_2E_1+E_2E_1^2=E_2^2E_1-(q^2+q^{-2})E_2E_1E_2+E_1E_2^2=0$.
The algebra $U_q^+({\mathfrak sl}_3)$ is usually named the
\emph{quantum Heisenberg algebra}. In the following we will denote
it by $\ha$. Setting $E_3=E_1E_2-q^2E_2E_1$, it is easy to check
(see for instance \cite{AD1}) that $\ha$ is the iterated Ore
extension generated over $\ku$ by the three generators
$E_1,E_2,E_3$ with relations:\bigbreak

\centerline{$E_1E_3=q^{-2}E_3E_1$, \qquad $E_2E_3=q^{2}E_3E_2 $,
\qquad $E_2E_1=q^{-2}E_1E_2 -q^{-2}E_3$.}\bigbreak

The center of $\ha$ is the a polynomial algebra in one
indeterminate $Z(\ha)=\ku[\Omega]$ where the quantum Casimir
element $\Omega$ is given by:\bigbreak

\centerline{$\Omega=(1-q^{-4})E_3E_1E_2+q^{-4}E_3^2=E_3\overline{E_3}$,
\quad with $\overline{E_3}=E_1E_2-q^{-2}E_2E_1$.}\bigbreak

Let $\ha=\bigoplus \ha_{m,n}$ be the canonical $\mathbb
N^2$-grading of $\ha$, defined putting $E_1$ on degree $(1,0)$ and
$E_2$ on degree $(0,1)$. In particular $E_3,\overline{E_3}\in
\ha_{1,1}$ and $\Omega\in \ha_{2,2}$.

\bigbreak
\subsubsection{Automorphisms and gradings of $U_q^+({\mathfrak{sl}_3})$}\label{autheissub}

For all $\alpha,\beta\in \ku^{\times}$, there exists one
automorphism $\widetilde{\psi}_{\alpha,\beta}$ of $\ha$ such that
$\widetilde{\psi}_{\alpha,\beta}(E_1)=\alpha E_1$ et
$\widetilde{\psi}_{\alpha,\beta}(E_2)=\beta E_2$. We introduce
$\widetilde{G}:=\{\widetilde{\psi}_{\alpha,\beta}\,;\,\alpha,\beta\in
\ku^{\times} \} \simeq(\ku^{\times})^2$, and the diagram
automorphism $\omega$ of $\ha$ defined by $\omega(E_1)=E_2$ and
$\omega(E_2)=E_1$. Studying the action of any algebra automorphism
of $\ha$ on the center and on the set of normal elements of $\ha$,
proposition 2.3 of \cite{AD1} prove that $\Aut_{\rm Alg} \ha$ is
the semi-direct
 product of $\widetilde{G}$ by the subgroup of order 2 generated by
$\omega$, (see  proposition 4.4 of \cite{C2} for another proof).
So we have for the type $A_2$ the following positive answer to
Problem \ref{pb}.

\begin{theo}\label{autA2}
For $\mathfrak g$ of type $A_2$, the algebra $\ha=U_q^+({\mathfrak
g})$ satisfies  $\Aut_{\rm Alg}\ha \simeq(\ku^{\times})^{2}
\rtimes \mathbb S_2$. \qed\end{theo}

Here is a consequence of this result which will be useful
in the next section.

\begin{corollary}\label{gradA2}
Let $\ha=\bigoplus \ha_{m,n}$ be the canonical $\N^2$-grading of
$\ha$. Let $\ha=\bigoplus T_{m,n}$ be another $\N^2$-algebra
grading of $\ha$. Then there exists a matrix
$\left(\begin{smallmatrix}p & q
\\ r&s \\ \end{smallmatrix}\right) \in M(2, \Z)$, with
non-negative entries, such that
$$\ha_{m,n}\subset T_{pm + qn, rm + sn}, \quad \text{ for all }(m,n)\in\N^2.$$
\end{corollary}

\pf This follows from Theorem \ref{fewt}, since $\ha$ has few
automorphisms by Theorem \ref{autA2}. \epf

\section{Partial results on the case where $\mathfrak g$ is of type $B_2$}
\label{b2}\medbreak A natural step in the study of Problem
\ref{pb} would be to consider the other Lie algebras $\mathfrak g$
of rank two. We summarize in this section some partial results
concerning the case $B_2$.

In this section, $\ku$ is an algebraically closed field of
characteristic zero, and $q\in\ku^{\times}$ a quantization
parameter not a root of one.

\bigbreak

\subsection{Some ring-theoretical properties of the algebra $U^+$}\label{rtp}

\subsubsection{Notations}\label{hyp}
 Let ${\mathfrak g}$ be the complex
simple Lie algebra over $\ku$ of type $B_2$. We denote
$U^+=U_q^+({\mathfrak g})$. We recall all notations of \ref{qea}
but denote now by $e_i$ the generators $E_i$; we have here $n=2$,
$C=\left(\begin{smallmatrix} 2  & -1  \\ -2  &  2  \\
\end{smallmatrix}\right)$, and $(d_1,d_2)=(2,1)$. Then the quantum
Serre relations are:
\medbreak

\centerline{ $\sum_{\nu=0}^{2} {\left[\begin{smallmatrix} 2 \\ \nu
\\ \end{smallmatrix}\right]}_{q^{2}} e_i^{2-\nu}e_j(-e_i)^{\nu}=0$
\quad for $i=1, j=2, a_{i,j}=-1, d_i=2$,} \bigbreak

\centerline{ $\sum_{\nu=0}^{3} {\left[\begin{smallmatrix} 3 \\ \nu
\\ \end{smallmatrix}\right]}_{q} e_i^{3-\nu}e_j(-e_i)^{\nu}=0$
\quad \ for $i=2, j=1, a_{i,j}=-2, d_i=1$.}
\medbreak

We compute the quantum binomial coefficients:
${\left[\begin{smallmatrix} 2 \\ 0  \\ \end{smallmatrix}\right]}_{q^{2}}=
{\left[\begin{smallmatrix} 2 \\ 2  \\ \end{smallmatrix}\right]}_{q^{2}}=1$,
${\left[\begin{smallmatrix} 2 \\ 1  \\ \end{smallmatrix}\right]}_{q^{2}}=
q^2+q^{-2}$,
${\left[\begin{smallmatrix} 3 \\ 0  \\ \end{smallmatrix}\right]}_{q}=
{\left[\begin{smallmatrix} 3 \\ 3  \\ \end{smallmatrix}\right]}_{q}=1$,
${\left[\begin{smallmatrix} 3 \\ 1  \\ \end{smallmatrix}\right]}_{q}=
{\left[\begin{smallmatrix} 3 \\ 2  \\ \end{smallmatrix}\right]}_{q}=
q^2+1+q^{-2}.$
We conclude that $U^+$ is the algebra generated over $\ku$ by two
generators $e_1$ and $e_2$ with commutation relations:
\medbreak

\centerline{$\left.\begin{matrix}
e_1^2e_2-(q^2+q^{-2})e_1e_2e_1+e_2e_1^2=0, \hfill &\quad{\rm (S1)}\hfill \\
e_2^3e_1-(q^2+1+q^{-2})e_2^2e_1e_2+
(q^2+1+q^{-2})e_2e_1e_2^2-e_1e_2^3=0.\hfill &\quad{\rm (S2)}\hfill \\
\end{matrix}\right.$}\bigbreak

\bigbreak
\subsubsection{$U^+$ as an iterated Ore extension}\label{ore}
From  the natural generators $e_1$ and $e_2$ of $U^+$, we
introduce following \cite{Ya} the $q$-brackets: \bigbreak

\centerline{$e_3=e_1e_2-q^2e_2e_1$ \qquad and \qquad
$z=e_2e_3-q^2e_3e_2$,} \bigbreak

Relations (S1) and (S2) imply: $e_1e_3=q^{-2}e_3e_1$, $e_1z=ze_1$
and $e_2z=ze_2$. In particular $z$ is central in $U^+$. From
\cite{Ya}, the monomials
$(z^ie_3^je_1^ke_2^l)_{(i,j,k,l)\in{\N}^4}$ form a PBW basis of
$U^+$. So $U^+$ is the algebra generated over $\ku$ by
$e_1,e_2,e_3,z$ with relations: \bigbreak

\centerline{
$\left.\begin{matrix}
e_3z=ze_3, \hfill& & \\
e_1z=ze_1, \hfill&\qquad e_1e_3=q^{-2}e_3e_1, \hfill& \\
e_2z=ze_2, \hfill&\qquad e_2e_3=q^{2}e_3e_2+z,
\hfill&\qquad e_2e_1=q^{-2}e_1e_2-q^{-2}e_3.\\
\end{matrix}\right.$}
\bigbreak

In other words $U^+$ is the iterated Ore extension ({\it cf.} \cite{BG}):
\medbreak

\centerline{
$U^+=\ku[e_3,z][e_1\,;\, \sigma][e_2\,;\, \tau,\delta]=S[e_2\,;\, \tau,\delta]$,
\quad with notation $S=\ku[e_3,z][e_1\,;\, \sigma]$,}\medbreak

\noindent where $\sigma$ is the automorphism of $\ku[e_3,z]$
defined by $\sigma(z)=z$, $\sigma(e_3)=q^{-2}e_3$, $\tau$ is the
automorphism of $\ku[e_3,z][e_1\,;\, \sigma]$ defined by
$\tau(z)=z$, $\tau(e_3)=q^{2}e_3$, $\tau(e_1)=q^{-2}e_1$, $\delta$
is the $\tau$-derivation of $\ku[e_3,z][e_1\,;\, \sigma]$ defined
by $\delta(z)=0$, $\delta(e_3)=z$, $\delta(e_1)=-q^{-2}e_3$, and
$S$ is the subalgebra of $U^+$ generated by $e_3$, $z$ and $e_1$.

\bigbreak
\subsubsection{Grading of $U^+$}\label{grad}
We consider the canonical grading $U^+=\bigoplus_{n\geq 0}U_n$
putting the natural generators $e_1$ and $e_2$ in degree one (and
then $e_3$ and $z$ are of degree 2 and 3 respectively) defined
from the basis $(z^ie_3^je_1^ke_2^l)_{(i,j,k,l)\in{\N}^4}$ de
$U^+$ by: ${U_n=\bigoplus_{3i+2j+k+l=n}\ku z^ie_3^je_1^ke_2^l}$
for any $n\geq 0$. We denote by $I={\bigoplus_{n\geq 1}U_n}$ the
ideal generated by $e_1,e_2,e_3,z$.

\bigbreak
\subsubsection{A localization of $U^+$}\label{loc}
The subalgebra of $U^+$ generated over $\ku$ by $e_1$ and $e_3$ is a quantum plane
(with $e_3e_1=q^2e_1e_3$); we will denote it by $\ku_{q^2}[e_3,e_1]$.
Its localization at the powers of $e_3$ and $e_1$ is the quantum torus
$\ku_{q^2}[e_3^{\pm},e_1^{\pm}]$.
The automorphism $\tau$ and the $\tau$-derivation $\delta$ extend to
$\ku_{q^2}[e_3^{\pm},e_1^{\pm}]$, and denoting by $V$ the algebra
$\ku_{q^2}[e_3^{\pm},e_1^{\pm}][z][e_2\,;\, \tau,\delta]$, we obtain the embedding:
\medbreak

\centerline{
$U^+=\ku[e_3,z][e_1\,;\, \sigma][e_2\,;\, \tau,\delta]=
\ku_{q^2}[e_3,e_1][z][e_2\,;\, \tau,\delta]
 \ \subset \
V=\ku_{q^2}[e_3^{\pm},e_1^{\pm}][z][e_2\,;\, \tau,\delta]$,}
\medbreak

Let us introduce in $U^+$ the bracket:
$w=e_2e_3-e_3e_2=z+(q^2-1)e_3e_2\in U_3$. It follows from
commutation relations in \ref{ore} that:
$e_1w=we_1+(1-q^{-2})e_3^2$, $e_2w=q^2we_2$, et $e_3w=q^{-2}we_3$.
Then the element:\bigbreak

\centerline{$z'=e_1w-q^{-4}we_1\in U_4$,}\bigbreak

satisfies $z'e_1=e_1z'$ and \ $z'e_2=e_2z'$, and so is central in $U^+$.
A straightforward computation shows that its development in the PBW basis of \ref{ore} is:
\medbreak

\centerline{$ z'=(1-q^{-4})(1-q^{-2})e_3e_1e_2
+q^{-4}(1-q^{-2})e_3^2+(1-q^{-4})ze_1$.}
\medbreak

In particular, $z'=s_1e_2+s_0$, with
$s_0=q^{-4}(1-q^{-2})e_3^2+(1-q^{-4})ze_1\in
\ku_{q^2}[e_3,e_1][z]$, and $s_1=(1-q^{-4})(q^{2}-1)e_1e_3$
non-zero in $\ku_{q^2}[e_3,e_1][z]$. So we have in $V$ the
identity $e_2=s_1^{-1}z'-s_1^{-1}s_0$, with $s_1^{-1}$ and
$s_1^{-1}s_0$ in $\ku_{q^2}[e_3^{\pm},e_1^{\pm}][z]$. Explicitly:
\medbreak

\centerline{$e_2=\frac1{(1-q^{-4})(q^2-1)}{e_3^{-1}e_1^{-1}z'}
+\frac1{q^4-1}{e_3^{-1}z}
-\frac1{q^2-1}e_1^{-1}e_3$\qquad\qquad\qquad(1)} \bigbreak

We conclude:

\centerline{$U^+\subset V=\ku_{q^2}[e_3^{\pm},e_1^{\pm}][z,z']$.}
\medbreak

Observe that $z$ and $z'$ being central in $V$, the only relation between the generators
of $V$ which is not a commutation is the $q^2$-commutation $e_3e_1=q^2e_1e_3$.

\bigbreak
\subsubsection{Conjugation in $U^+$}\label{conj}
We have introduced in \ref{ore} the homogeneous element $e_3$ of
degree two defined from natural generators $e_1$ et $e_2$ by
$e_3=e_1e_2-q^2e_2e_1$, which satisfy $e_1e_3=q^{-2}e_3e_1$  and
$e_2e_3-q^2e_3e_2=z$. Conjugating $q$ in $q^{-1}$, we can also
consider:

$$\overline{e_3}=e_1e_2-q^{-2}e_2e_1=(1-q^{-4})e_1e_2+q^{-4}e_3,$$

\noindent and prove that
$e_1\overline{e_3}=q^{2}\overline{e_3}e_1$  and
$e_2\overline{e_3}-q^{-2}\overline{e_3}e_2=q^{-4}z$. We obtain in
particular the relations
$e_3\overline{e_3}=(1-q^{-4})q^2e_1e_3e_2+q^{-4}e_3^2$ and
$z'=(1-q^{-2})(e_3\overline{e_3}+(1+q^{-2})ze_1)$, which will be
used further in the paper.

\bigbreak
\subsubsection{Center and normalizing elements of $U^+$}\label{centnorm}
Denote by $Z(U^+)$ the center and $N(U^+)$ the set of normalizing elements in  $U^+$

\medbreak
\begin{lemma}\label{cent} $N(U^+)=Z(U^+)=\ku[z,z']$.\end{lemma}
\medbreak

\pf The calculation of $Z(U^+)$ can be deduced from general
results of \cite{C1}. The equality $N(U^+)=Z(U^+)$ for the type
$B_2$ was observed at \cite{C2}, remark 2.2 (iii). We give here a
short direct proof using the embedding $U^+\subset V$. Take $f\in
N(U^+)$ non-zero. From proposition 2.1 de \cite{C2}, $f$ is
$q$-central in $U^+$, that is there exist $m,n\in{\mathbb Z}$ such
that $fe_1=q^me_1f$ and $fe_2=q^ne_2f$, and so $fe_3=q^{m+n}e_3f$.
In $V$, the element $f$ is a finite sum: \medbreak

\centerline{ $f=\sum_{i,j\in{\Z}}f_{i,j}(z,z')e_1^ie_3^j$ \
\ with $f_{i,j}(z,z')\in \ku[z,z']$.}
\medbreak

Since the polynomials $f_{i,j}(z,z')$ are central in $V$, the
identities $fe_1=q^me_1f$ and $fe_3=q^{m+n}e_3f$ give by
identification $2j=m$ and $2i=-m-n$ for all $i,j$ such that
$f_{i,j}(z,z')\not=0$. Then $f=f_{i,j}(z,z')e_1^ie_3^j$, where
$i=-\frac{m+n}{2}$ and $j=\frac{m}{2}$. Its follows from relation
(1) of  \ref{loc} and from $q$-commutation $fe_2=q^ne_2f$ that
$-2j+2i=2i=-2j-2i=n$. So $i=j=n=0$, and $f=f_{0,0}(z,z')\in
\ku[z,z']$. We have proved that $N(U^+)\subseteq \ku[z,z']$. The
inverse inclusions $\ku[z,z']\subseteq Z(U^+)\subseteq N(U^+)$ are
clear.\epf

\begin{obs}\label{rks} For $\mathfrak
g$ of type $B_2$,  $\Autdiag ({\mathfrak g})$ is trivial and so
the main Problem \ref{pb} must be here formulated as: do we have
$\Aut_{\rm Alg} U^+\simeq(\ku^{\times})^{2}$ ? The method used in
\cite{AD1} for the case $A_2$ was based on the facts that any
automorphism of $\ha$ preserves the center $Z(\ha)$ which is a
polynomial algebra in one variable and the non-empty set of normal
but non central elements of $\ha$. It follows from lemma
\ref{cent} that the second argument fails for the case $B_2$, and
that the first one is much more complicated to use in view of the
structure of automorphism group of a commutative polynomial
algebra in two variables. A natural idea to determine the group
$\Aut_{\rm alg}U^+$ is to study the action of an automorphism on
the prime spectrum of $U^+$. The structure of this spectrum
remains widely unknown as far as we know (see further final remark
of \ref{z'}) and we only present some partial results in the
following. In particular the two central generators $z$ and $z'$
do not play symmetric r\^oles (see proposition \ref{zz'}) and we
conjecture that any automorphism of $U^+$ stabilizes the prime
ideal $(z)$, which would be sufficient to solve the problem
because of proposition \ref{autstabz}.
\end{obs}

\subsection{Automorphisms stabilizing the prime ideal
$(z)$}\label{z}

\subsubsection{The factor algebra $U^+/(z)$}\label{factalg}
It is clear that the ideal $(z)$ generated in $U^+$ by the central
element $z$ is completely prime, and the factor domain $U^+/(z)$
is  the quantum enveloping algebra $\ha=U_q^+({\mathfrak sl}_3)$
considered in \ref{heis}. The canonical map $\pi : U^+\rightarrow
U^+/(z)$ induces an isomorphism between $U^+/(z)$ and $\ha$
defined by $\pi(e_1)=E_1$ and $\pi(e_2)=E_2$, where $E_1$ and
$E_2$ are the natural generators of $\ha$ introduced in
\ref{heis}. Observe that $\pi(e_3)=E_3$ and
$\pi(z')=(1-q^{-2})\Omega$.

\bigbreak
\subsubsection{The group $\Aut_z(U^+)$}\label{autz}
We introduce the subgroup $\Aut_z(U^+)$ of all automorphisms of
the algebra $U^+$ which stabilize the ideal $(z)$. In particular
$\Aut_z(U^+)$ contains the subgroup
$G:=\{\psi_{\alpha,\beta}\,;\,\alpha,\beta\in \ku^{\times}
\}\simeq(\ku^{\times})^2$ where $\psi_{\alpha,\beta}$ is the
automorphism defined for any $\alpha,\beta\in \ku^{\times}$ by:
\bigbreak

\centerline{$\psi_{\alpha,\beta}(e_1)=\alpha e_1$ \quad et \quad
$\psi_{\alpha,\beta}(e_2)=\beta e_2$.}\medbreak

It is clear that any $\theta\in \Aut_z U^+$ induces an
automorphism $\widetilde{\theta} \in \Aut \ha$ defined by
$\widetilde{\theta}(E_1)=\pi(\theta(e_1))$ et
$\widetilde{\theta}(E_2)=\pi(\theta(e_2))$. We denote by $\Phi$
the group morphism $\Phi: \Aut U^+\rightarrow \Aut \ha \ ;\ \theta
\mapsto \widetilde{\theta}$. The notation is coherent since
$\Phi(G)=\widetilde{G}$ for the group $\widetilde
G\simeq(\ku^{\times})^2$ introduced in \ref{heis}. The next
proposition proves that $\Phi$ defines an isomorphism between
$\Aut_zU^+$ and $\widetilde G$.\medbreak

\begin{proposition}\label{autstabz}
The subgroup $\Aut_z(U^+)$ of algebra automorphisms of $U^+$
stabilizing the ideal  $(z)$ is isomorphic to $(\ku^{\times})^2$.
\end{proposition}

\pf {\it Step 1:} we prove that, for any $\theta\in\Aut_z U^+$,
there exist $\lambda,\mu\in \ku^{\times}$ and $p(z)\in \ku[z]$
such that $\theta(z)=\lambda z$ and $\theta(z')=\mu
z'+p(z)$.\bigbreak

Take $\theta\in\Aut_z U^+$. There exists $u\in U^+$, $u\not=0$
such that $\theta(z)=uz$. But $\theta(z)\in Z(U^+)$ because $z\in
Z(U^+)$, and then $u\in Z(U^+)$. Similarly $\theta^{-1}(z)=vz$ for
some $v\in Z(U^+)$, $v\not=0$. So
$z=\theta(\theta^{-1}(z))=\theta(v)uz$ in $Z(U^+)=\ku[z,z']$ (see
\ref{centnorm}), which implies $u\in \ku^{\times}$. Denoting $u=\lambda$,
we conclude that $\theta(z)=\lambda z$ with $\lambda\in \ku^{\times}$.
The restriction of $\theta$ to $Z(U^+)$ is a $\ku$-automorphism of
$\ku[z,z']$ such that $\theta(z)=\lambda z$. By surjectivity, the
$z'$-degree of $\theta(z')$ is necessarily 1. Denote
$\theta(z')=r(z)z'+p(z)$ with $r(z),p(z)\in \ku[z]$, $r(z)\not=0$.
Using an analogue expression for $\theta^{-1}(z')$ in the equality
$z'=\theta(\theta^{-1}(z'))$ we obtain that $r(z)=\mu\in\ku^{\times}$.
\medbreak

{\it Step 2:} we prove that $\Ima\Phi=\widetilde{G}$. \bigbreak

Let $\theta\in \Aut U^+$ and $\widetilde{\theta}=\Phi(\theta)$.
From  Theorem \ref{autA2}, there exist $\alpha,\beta\in
\ku^{\times}$ and $i\in\{0,1\}$ such that
$\widetilde{\theta}=\widetilde{\psi}_{\alpha,\beta}\omega^{i}$.
Suppose that $i=1$. Then $\omega=\Phi(\theta')$ for
$\theta'={\psi}_{\alpha^{-1},\beta^{-1}}\theta$, which satisfies
$\pi(\theta'(e_1))=\omega({\text{\sc e}}_1)={\text{\sc e}}_2$ and
$\pi(\theta'(e_2))=\omega({\text{\sc e}}_2)={\text{\sc e}}_1$. In
other words, there exist $a,b\in U^+$ such that
$\theta'(e_1)=e_2+za$ and $\theta'(e_2)=e_1+zb$. Applying
$\theta'$ to the first Serre relation (S1) in $U^+$, we obtain:
\bigbreak

\centerline{
$(e_2+za)^2(e_1+bz)-(q^2+q^{-2})(e_2+za)(e_1+bz)(e_2+za)+(e_1+bz)(e_2+za)^2=0$.}
\bigbreak

Using the grading of \ref{grad}, this identity develops into an
expression  $s+t=0$ with
$s=e_2^2e_1-(q^2+q^{-2})e_2e_1e_2+e_1e_2^2\in U_3$ and the rest
$t$ in $\bigoplus_{n\geq 5}U_n$. Then $s=0$. But the relations of
\ref{conj} allow to compute:
$s=e_2(e_2e_1-q^2e_1e_2)-q^{-2}(e_2e_1-q^2e_1e_2)e_2 =
-q^2e_2\overline{e_3}+\overline{e_3}e_2 = -q^{-2}z\not=0$. So a
contradiction. We conclude that $i=0$, and then
$\widetilde{\theta}=\widetilde{\psi}_{\alpha,\beta}\in\widetilde
G$.\medbreak

{\it Step 3:} we prove that $\Phi$ is injective.\bigbreak

We fix an automorphism $\theta\in\Aut_z(U^+)$ such that
$\theta\in\Ker\Phi$. By definition of $\Ker\Phi$, there exist
$a,b\in U^+$ such that: \bigbreak

\centerline{ $\theta(e_1)=e_1+za$\quad and \quad
$\theta(e_2)=e_2+zb$.}\bigbreak

The elements $e_3$ and $\overline{e_3}$ being defined as
$q$-brackets of $e_1$ and $e_2$, we deduce: \bigbreak

\centerline{$\begin{matrix} \theta(e_3)=e_3+zc, &{\rm where}
\ c=(e_1b-q^2be_1)+(ae_2-q^2e_2a)+z(ab-q^2ba)\hfill\\
\theta(\overline{e_3})=\overline{e_3}+z\overline{c}, &{\rm where}
\ \overline{c}=(e_1b-q^{-2}be_1)+(ae_2-q^{-2}e_2a)+z(ab-q^{-2}ba)\hfill\\
\end{matrix}$}\bigbreak

Applying $\theta$ to the relations $e_1e_3=q^{-2}e_3e_1$ (see
\ref{ore}) and $\overline{e_3}=(1-q^{-4})e_1e_2+q^{-4}e_3$ (see
\ref{conj}) we obtain by identification:\bigbreak

\centerline{$\begin{matrix}(e_1c-q^{-2}ce_1)+(ae_3-q^{-2}e_3a)+z(ac-q^{-2}ca)=0 \hfill\\
(1-q^{-4})ae_2+(1-q^{-4})e_1b+q^{-4}c+(1-q^{-4})zab=\overline{c}\hfill\\
\end{matrix}$}\bigbreak

Consider now $z=e_2e_3-q^2e_3e_2$. By step one, there exists some
$\lambda\in \ku^{\times}$ such that: $\lambda
z=z+z(e_2c-q^2ce_2)+z(be_3-q^2e_3b)+z^2(bc-q^2cb)$. Simplifying by
$z$, we obtain
$\lambda-1=(e_2c-q^2ce_2)+(be_3-q^2e_3b)+z(bc-q^2cb)$. The right
member lies in the ideal $I$ defined in \ref{grad} and the left
member is a scalar in $\ku=U_0$. Therefore $\lambda=1$. We
conclude that:\bigbreak

\centerline{$\theta(z)=z$ \quad and \quad
$(e_2c-q^2ce_2)+(be_3-q^2e_3b)+z(bc-q^2cb)=0$.}\bigbreak

Similar calculations for the other central generator
$z'=(1-q^{-2})(e_3\overline{e_3}+(1+q^{-2})ze_1)$ (see
\ref{conj}), using the equality $\theta(z')=\mu z'+p(z)$ from step
one, give $\mu=1$ and $p(z)=(1-q^{-2})zs(z)$ for some
$s(z)\in\ku[z]$ which satisfies:\bigbreak

\centerline{
$s(z)=(1+q^{-2})za+zc\overline{c}+e_3\overline{c}+q^{-4}ce_3+
(1-q^{-4})ce_1e_2$.} \bigbreak

From  \ref{ore}, we can consider the degree function $\deg$ at the
indeterminate $e_2$ in the polynomial algebra $U^+=S[e_2\,;\,
\tau,\delta]$. Introduce in particular the three positive integers
$d=\deg\theta(e_1)$, $d'=\deg\theta(e_2)$ and
$d''=\deg\theta(e_3)$. Comparing the degree of the two members of
all the equalities obtained above, we obtain by very technical
considerations whose details are left to the reader that we have
necessarily $d=d''=0$ et $d'=1$. In other words there exist
$a,b_1,b_0,c\in S$ such that:

\centerline{ $\theta(e_1)=e_1+za$,\quad
$\theta(e_2)=e_2+z(b_1e_2+b_0)$,\quad $\theta(e_3)=e_3+zc$, \quad
$\theta(z)=z$.} \bigbreak

In particular the restriction $\theta_S$ of $\theta$ to $S$ is an
automorphism of $S$ fixing $z$. Consider the field of fractions
$K=\ku(z)$, the algebra $T=K[e_3][e_1\,;\,\sigma]\supseteq S$ and
the extension $\theta_T$ of $\theta_S$ to $T$. Since $T$ is a
quantum plane over $\ku$ (with $e_1e_3=q^{-2}e_3e_1$), we can
apply proposition 1.4.4 of\cite{AC} and deduce that the
$K$-automorphism $\theta_T$ satisfies $\theta(e_1)=fe_1$ and
$\theta(e_3)=ge_3$ for some $f,g\in \ku^{\times}$. Because
$\theta(e_1)$ and $\theta(e_3)$ are in the subalgebra
$S=\ku[z][e_3][e_1\,;\,\sigma]$ of
$T=\ku(z)[e_3][e_1\,;\,\sigma]$, we have in fact $f,g$ non-zero in
$\ku[z]$. By the same argument for the automorphism $\theta^{-1}$,
there exist $f',g'$ non-zero in $\ku[z]$ such that
$\theta^{-1}(e_1)=f'e_1$ and $\theta^{-1}(e_3)=g'e_3$. By
composition of $\theta$ and $\theta^{-1}$ it follows that
$ff'=gg'=1$, and then $f,g\in \ku^{\times}$. Denoting $f=\alpha$
and $g=\gamma$, we obtain $\alpha e_1=e_1+za$ and $\gamma
e_3=e_3+zc$. These equalities in $S$ imply  $(\alpha-1)e_1\in zS$
and $(\gamma-1)e_3\in zS$ with $\alpha,\gamma\in \ku^{\times}$,
and so $\alpha=\gamma=1$ and $a=c=0$. To sum up, we have
$\theta(e_1)=e_1$, $\theta(e_3)=e_3$ and $\theta(z)=z$. It is then
easy to check that we have also $\theta(e_2)=e_2$. We conclude
that $\theta=\id_{U^+}$.\epf \medbreak

\subsection{Non permutability of the central generators $z$ and $z'$}\label{z'}

\subsubsection{The prime ideal \,$(z')$}
The ideal $(z')$ generated in $U^+$ by the central element $z'$ is
completely prime. A direct proof consists in checking by
computations in the iterated Ore extension $U^+$ and its
localization $V=\ku_{q^2}[e_3^{\pm 1},e_1^{\pm 1}][z,z']$ that any
element $v\in V$ satisfying $z'v\in U^+$ is necessarily an element
of $U^+$; then the ideal $(z')=z'U^+$ is no more than the
contraction $z'V\cap U^+$ of the completely prime ideal $z'V$ of
$V$. The complete primeness of $(z')$ can also be proved by the
algorithmic method of \cite{Cau}, or can be deduced from the
description of the $\mathcal H$-prime spectrum of $U^+$ (see
further \ref{hspec}) following the method of \cite{G}. \medbreak

\begin{proposition}\label{zz'}
There are no algebra automorphisms of $U^+$ sending $(z)$ to
$(z')$ or $(z')$ to $(z)$.
\end{proposition}\medbreak

\pf Recall that $\ha$ denotes the factor algebra $U^+/(z)$ and
$\pi$ the canonical map $U^+\rightarrow \ha$ (see \ref{heis} and
\ref{factalg}). Denote $\ha'=U^+/(z')$ and $\pi'$ the canonical
map $U^+\rightarrow \ha'$. We have observed that the center of
$\ha$ is the polynomial algebra $Z(\ha)=\ku[\Omega]$; note that
$\Omega$ equals $\pi(z')$ up to a multiplication by a non-zero
scalar. By direct calculations in $\ha'$ one can prove similarly
that the center of $\ha'$ is the polynomial algebra
$Z(\ha')=\ku[\pi'(z)]$. Consider the standard $\N^2$-grading
$U^+=\bigoplus U_{m,n}$ putting $e_1$ on degree $(1,0)$ and $e_2$
on degree $(0,1)$. The central generator $z$ is homogeneous of
degree $(1,2)$ and the $\N^2$-grading induced on $\ha=U^+/(z)$ is
just the $\N^2$-grading $\ha=\bigoplus \ha_{n,m}$ considered in
\ref{heis}. The central generator $z'$ is homogeneous of degree
$(2,2)$ and the $\N^2$-grading induced on $\ha'=U^+/(z')$ will be
denoted by $\ha'=\bigoplus \ha'_{n,m}$.\bigbreak

Now suppose that there exists some algebra automorphism $\theta$
of $U^+$ such that $\theta(z')=(z)$. Then $\theta$ induces an
isomorphism $\widehat\theta:\ha'\rightarrow \ha$ defined by
$\widehat\theta\pi'=\pi\theta$. In particular $\ha$ can be graded
by the $\N^2$-grading $\ha=\bigoplus T_{n,m}$ defined by
$T_{n,m}=\widehat\theta(\ha'_{m,n})$. Since $E_3,\overline{E_3}\in
\ha_{1,1}$ and $\Omega\in \ha_{2,2}$, it follows from corollary
\ref{gradA2} there exist two integers $r,s\geq 1$ such that
$E_3,\overline{E_3}\in T_{r,s}$ and $\Omega\in T_{2r,2s}$. Set
$t_1=\widehat\theta^{-1}(E_3)$,
$t_2=\widehat\theta^{-1}(\overline{E_3})$, and
$t=\widehat\theta^{-1}(\Omega)=\widehat\theta^{-1}(E_3\overline{E_3})=t_1t_2$.
By construction we have $t_1,t_2\in \ha'_{r,s}$ and $t\in
\ha'_{2r,2s}$. Because $Z(\ha)=\ku[\Omega]$ and $\widehat\theta$
is an isomorphism of algebras, $t$ is a generator of the
polynomial algebra $Z(\ha')=\ku[\pi'(z)]$. Then there exist
$\lambda\in\ku^{\times}, \mu\in\ku$ such that $t=\lambda
\pi'(z)+\mu$. Since $t\in \ha'_{2r,2s}$, we have necessarily
$\mu=0$. We obtain in $\ha'$ the equality
$\pi'(z)=\lambda^{-1}t_1t_2$ where $\lambda\in\ku^{\times}$ and
$t_1,t_2\in \ha'_{r,s}$. Choosing $u_1,u_2\in U_{r,s}$ such that
$\pi'(u_1)=\lambda^{-1}t_1$ and $\pi'(u_2)=t_2$, we deduce in
$U^+$ an equality $z=u_1u_2+z'u$ for some $u\in U^+$. Back to the
$\N$-grading $U^+=\bigoplus U_n$ introduced in \ref{grad}, we have
clearly $u_1u_2\in U_{2(r+s)}$ whereas $z\in U_3$ and $z'\in U_4$.
It follows that $u=0$, and then $z=u_1u_2$ gives the
contradiction.\epf

\begin{obs} Suppose that we can
prove similarly that $U^+/I$ is not isomorphic to $\ha\simeq
U^+/(z)$ for any height one prime ideal $I$ of $U^+$, then we
could deduce that any algebra automorphism of $U^+$ necessarily
stabilizes $(z)$ and then give by proposition \ref{autstabz} a
positive answer to Problem \ref{pb} (see \ref{rks}). For example,
if $I=(z-\alpha)$ with $\alpha\in\ku^{\times}$, it is possible to
separate $U^+/I$ from $\ha$ up to isomorphism (by technical
considerations on the $q$-brackets in the both factor algebras
which are not developed here). Unfortunately the complete
description of height one prime ideals in $U^+$ is not known as
far as we know. We restrict in the last subsection to the graded
prime ideals.

\end{obs}

\subsection{Action on the $\mathcal H$-spectrum}\label{hspec}

\subsubsection{$\mathcal H$-spectrum of $U_q^+(\mathfrak g)$}\label{gor}
In this subsection, $\mathfrak g$ is an arbitrary simple
finite-dimensional Lie algebra.
 We consider all data and
notations of \ref{qea}.
 Denote by $\mathcal H$ the $n$-dimensional torus
$(\ku^{\times})^n$. The canonical $\Z^n$-grading on
$U_q^+(\mathfrak g)$, given by $\deg E_i=\epsilon_i$, the $i$-th.
term of the canonical basis of $\Z^n$, induces a rational action
of $\mathcal H$ on $U_q^+(\mathfrak g)$. The $\mathcal H$-spectrum
of $U_q^+(\mathfrak g)$ is the set of graded prime ideals. We
denote it by $\mathcal H$-$\spec U_q^+(\mathfrak g)$. It
determines a stratification of the whole prime spectrum $\spec
U_q^+(\mathfrak g)$. We refer to \cite[section II]{BG} for more
details on stratifications of iterated Ore extensions, and
summarize here some results obtained in \cite{G} concerning the
case of $U_q^+(\mathfrak g)$.

\bigbreak More generally, \cite{G} determines a stratification of
$\spec S^{w}$ for a family of algebras $S^{w}$ (denoted also by
$R_0^{w}$, see  \cite[p. 217, l. 7]{G} for the equivalence between
the two notations), where $w$ is an element of the Weyl group $W$
of $\mathfrak g$.  The space of graded prime ideals of $S^{w}$ is
indexed by the set
$$W {\diamondsuit}^{\omega} W := \{(x,y) \in
W\times W: x \le w \le y\},$$ where $\le$ is the Bruhat order (see
for example \cite[App. 1]{J}). The algebra $U_q^+(\mathfrak g)$ is
the particular case $w=e$, and thus the index set of $\mathcal
H$-$\spec U_q^+(\mathfrak g)$ is just $W$. The author introduces
in \cite{G}  an ideal $Q(y)$ for any $y\in W$ (this is the ideal
$Q(e,y)_e$ with the notations of \cite[p. 231 and p. 236]{G}), and
proves (see \cite[7.1.2, (ii), p. 239]{G} that
$$\mathcal H-\spec U_q^+(\mathfrak g)=\{Q(y)\}_{y\in
W}.$$ By \cite[6.11, p. 236]{G}, we have $Q(y) \subset Q(y')$ if
and only if $y' \le y$. Furthermore, \cite[5.3.3, p. 225]{G} (with
notation $Y(y)=Y_e(e,y)$) asserts that
\begin{equation}\label{stratif}
\spec U_q^+(\mathfrak g) = \coprod_{y\in W} Y(y),
\end{equation}
\noindent where $Y(y)$ denotes the set of prime ideals containing
$Q(y)$. By \cite[6.12, p. 237]{G}, each subset $Y(y)$ has a unique
minimal element, namely $Q(y)$, and coincides with the $\mathcal
H$-strata of $\spec U_q^+(\mathfrak g)$ corresponding to the
graded prime ideal $Q(y)$, that is:
$$Y(y)=\spec_{Q(y)}U_q^+(\mathfrak g)=\{P\in\spec U_q^+(\mathfrak
g)\,:\, Q(y)=\bigcap_{h\in\mathcal H}h.P\}.$$ Applying \cite[6.13,
p. 238]{G}, the closure of each $Y(y)$ is a union of other
$Y(y')$'s, namely $\overline{Y(y)} = \coprod_{y'\in W, y'\le y}
Y(y'),$ and the disjoint union \eqref{stratif} is then a
stratification of $\spec U_q^+(\mathfrak g)$.

\subsubsection{$\mathcal H$-spectrum of $U^+$ for the type
$B_2$}\label{hspecb2} Suppose now that $\mathfrak g$ is of type
$B_2$ and keep all notations of \ref{rtp}. The Weyl group $W$ is
of order 8. Its elements  can be described by their action on the
roots $\{\varepsilon_1,\varepsilon_2\}$ and as words on the two
generators $s_1$ and $s_2$ as follows:
$$\begin{matrix}
\hfill s_1: &\varepsilon_1 \mapsto \varepsilon_2,\hfill
 &\varepsilon_2 \mapsto\varepsilon_1,\hfill &\qquad\qquad
 \hfill s_2: &\varepsilon_1 \mapsto \varepsilon_1,\hfill
  &\varepsilon_2 \mapsto-\varepsilon_2,\hfill \\
\hfill s_1s_2:&\varepsilon_1 \mapsto \varepsilon_2,\hfill
&\varepsilon_2 \mapsto-\varepsilon_1,\hfill &\qquad\qquad
\hfill s_2s_1: &\varepsilon_1 \mapsto -\varepsilon_2,\hfill
  &\varepsilon_2 \mapsto\varepsilon_1,\hfill \\
\hfill s_1s_2s_1s_2=s_2s_1s_2s_1: &\varepsilon_1 \mapsto -\varepsilon_1,
\hfill  &\varepsilon_2 \mapsto\varepsilon_2,\hfill &\qquad\qquad
 \hfill e: &\varepsilon_1 \mapsto \varepsilon_1,\hfill
  &\varepsilon_2 \mapsto\varepsilon_2.\hfill \\
\end{matrix}$$

Using the results of \cite{G} recalled in \ref{gor}, the $\mathcal
H$-prime spectrum of $U^+=U_q^+(\mathfrak g)$ has exactly 8
elements.
 The ideals $(0)$, $(e_1)$, $(e_2)$ and $(e_1,e_2)$ are clearly
graded prime ideals of $U^+$. This is also the case for the ideals
$(e_3)$ and $(\overline{e_3})$, the factor algebras being in both
cases domains isomorphic to a quantum plane. Finally the prime
ideals $(z)$ and $(z')$ considered in \ref{z} and \ref{z'} are
also graded prime ideals of $U^+$. So the poset $(W,\leq)$ and the
$\mathcal H$-spectrum of $U^+$ are:

$$
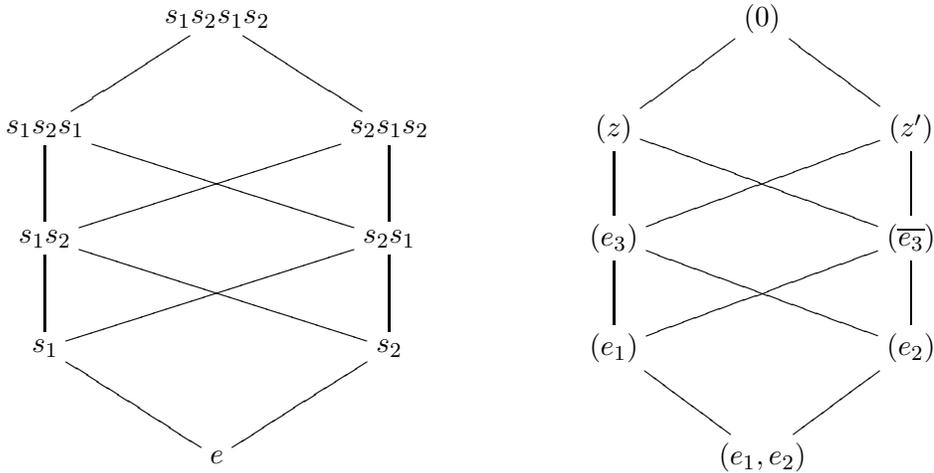
\xymatrix{
  & s_1s_2s_1s_2\ar@{-}[dl]\ar@{-}[dr] & & & &(0)\ar@{-}[dl]\ar@{-}[dr] & \\
s_1s_2s_1\ar@{-}[d]\ar@{-}[drr] &  &s_2s_1s_2 \ar@{-}[d]\ar@{-}[dll] & & (z)\ar@{-}[d]\ar@{-}[drr] & &(z')\ar@{-}[d]\ar@{-}[dll]\\
s_1s_2 \ar@{-}[d]\ar@{-}[drr] & & s_2s_1\ar@{-}[d]\ar@{-}[dll]
& & (e_3)\ar@{-}[d]\ar@{-}[drr] & &(\overline{e_3})\ar@{-}[d]\ar@{-}[dll]\\
s_1\ar@{-}[dr]  & & s_2\ar@{-}[dl] & & (e_1)\ar@{-}[dr]  & &(e_2)\ar@{-}[dl]\\
&  e & & &  &(e_1,e_2) & \\}$$

\bigbreak
\begin{proposition}\label{authspec}
The subgroup of algebra automorphisms of $U^+$ stabilizing the
$\mathcal H$-spectrum of $U^+$ is isomorphic to $(\ku^{\times})^2$.
\end{proposition}

\pf Any automorphism of $U^+$ stabilizing the $\mathcal
H$-spectrum of $U^+$ preserves the set of height one prime graded
ideals of $U^+$, that is $\{(z),(z')\}$. Then the result follows
from Propositions \ref{autstabz} and \ref{zz'}.\epf

\medbreak
\subsubsection{Automorphisms of some subalgebras of $U^+$ indexed by
$W$.}\label{subalg} Take all data and notations of \ref{qea} and
consider for any $w$ in the Weyl group $W$ of $\mathfrak g$ and
for $e=\pm 1$ the subalgebras $U^+(w,e)$ of $U_q^+(\mathfrak g)$
defined in \cite{L}. For any reduced expression
$w=s_{i_1}s_{i_2}\ldots s_{i_n}$ of an element $w\in W$ the
elements
$$E_{i_1}^{(c_1)}T'_{i_1,e}(E_{i_2}^{(c_2)})\cdots
T'_{i_1,e}T'_{i_2,e}\cdots T'_{i_{n-1},e}(E_{i_n}^{(c_n)})$$ for
various $(c_1,c_2,\ldots,c_n)\in{\bf N}^n$ form a basis of a
subspace $U^+(w,e)$ of $U_q^+(\mathfrak g)$ which does not depend
of the reduced expression of $w$ (\cite{L} 40.2.1 p. 321). The
Lusztig automorphisms $T'_{i,e}$ appearing in this definition are
the symmetries of $U_q^+(\mathfrak g)$ related to the braid group
action defined in \cite{L} 37.1.2. From \cite{L} 40.2.1 (d), we
have: $\ell(s_iw)=\ell(w)-1\Rightarrow E_iU^+(w,e)\subset
U^+(w,e)$. In particular, for $w_0$ the element of maximal length
in $W$, we have $U^+(w_0,e)=U_q^+(\mathfrak g)$ (see the remark
after 40.2.2 of \cite{L}). The subspaces $U^+(w,e)$ can be
identified with the subalgebras $U_q(n_w)$ of lemma 1.5 of
\cite{C2} and also appear in \cite{J} p. 123.\bigbreak

Suppose now that $\mathfrak g$ is of type $B_2$, take all
notations of \ref{rtp} and denote by $A_w$ the subalgebra
$U^+(w,1)$ of $U^+$, for any $w\in W$ expressed as a word into the
generators $s_1$ and $s_2$ of \ref{hspecb2}. Using the above
definition of the basis of $A_w$ and the properties of the
automorphisms $T'_{j,1}$ (in particular \cite{L} 39.2.3),
straightforward calculations allow to obtain the following
description by generators of the eight subalgebras $A_w$ of
$U^+$:\bigbreak

$$\xymatrix{
  &{\begin{matrix}
U^+=A_{s_1s_2s_1s_2}\\
\quad\ \ =A_{s_2s_1s_2s_1}\\
\end{matrix}}
 \ar@{-}[dl]\ar@{-}[dr] & \\
{\ku}\langle e_2,w,e_3\rangle=A_{s_1s_2s_1}\ar@{-}[d] &
& A_{s_2s_1s_2}={\ku}\langle e_1,\overline{e_3},\overline{w}\rangle\ar@{-}[d]\\
{\ku}\langle e_2,w\rangle=A_{s_1s_2} \ar@{-}[d] & & A_{s_2s_1}={\ku}\langle e_1,\overline{e_3}\rangle\ar@{-}[d] \\
{\ku}\langle e_2\rangle=A_{s_1}\ar@{-}[dr]  & & A_{s_2}={\ku}\langle e_1\rangle\ar@{-}[dl] \\
& A_e=\ku & \\}$$

\bigbreak
with the commutation relations:\medbreak

$\begin{matrix} {\rm left\ side}\hfill &  {\left\{\begin{matrix}
e_2w=q^2we_2, \hfill& & \\
e_3w=q^{-2}we_3, \hfill&\qquad e_3e_2=e_2e_3-w, \hfill& \\
e_1w=we_1+(1-q^{-2})e_3^2, \hfill&\qquad e_1e_2=q^{2}e_2e_1+e_3,
\hfill&\qquad e_1e_3=q^{-2}e_3e_1.\\
\end{matrix}\right.}\hfill\\
{ } & \\
 {\rm right\ side}\hfill & {\left\{\begin{matrix}
e_1\overline{e_3}=q^2\overline{e_3}e_1, \hfill& & \\
\overline{w}\,\overline{e_3}=q^{-2}\overline{e_3}\,\overline{w},
\hfill&\qquad
\overline{w}e_1=e_1\overline{w}+(q^2-1)\overline{e_3}^2, \hfill& \\
e_2\overline{e_3}=\overline{e_3}e_2+\overline{w}, \hfill&\qquad
e_2e_1=q^2e_1e_2-q^2\overline{e_3},
\hfill&\qquad e_2\overline{w}=q^{-2}\overline{w}e_2\\
\end{matrix}\right.}\hfill\\
\end{matrix}$
\bigbreak \medbreak

The eight algebras are iterated Ore extensions. At level one,
$A_{s_1}$ and $A_{s_2}$ are just commutative polynomial algebras
in one variable. At level two, $A_{s_1s_2}$ and $A_{s_2s_1}$ are
isomorphic to a same quantum plane (with parameter $q^2$) and so
by \cite{AC} their group of algebra automorphisms is isomorphic to
$(\ku^{\times})^2$. The third level introduces some asymmetry in the
diagram. One can prove by direct calculations (which are left to
the reader) similar to the proof of \cite[Lemme 2.2 and Proposition 2.3]{AD1} that:
 \begin{itemize}
 \item[(i)] the center of the algebra $A_{s_1s_2s_1}$ is $Z=\ku[z]$
 with $z=(1-q^2)e_3e_2+w$;\bigbreak

 \noindent its set of normal elements is
 $N=\bigcup_{n\geq 0}\ku[z]w^n$;\bigbreak

 \noindent its automorphism group is $(\ku^{\times})^2$ acting by $(\alpha,\beta):e_2\mapsto
 \alpha e_2, \ e_3\mapsto \beta e_3,\ w\mapsto \alpha\beta w$.\medbreak

 \item[(ii)] the center of the algebra $A_{s_2s_1s_2}$ is $Z'=\ku[u]$
 with $u=(1-q^{-4})e_1\overline{w}+(q^2-1)\overline{e_3}^2$;\medbreak

\noindent  its set of normal elements is
 $N'=\bigcup_{n\geq }\ku[u]\overline{e_3}^n$;\medbreak

\noindent  its automorphisms group
 is  $(\ku^{\times})^2$ acting by $(\alpha,\gamma):e_1\mapsto
 \alpha e_1,\ \overline{e_3}\mapsto \gamma \overline{e_3},\
  \overline{w}\mapsto \alpha^{-1}\gamma^2\overline{w}$.
\end{itemize}

\bigbreak Using the fact that an isomorphism from
$A_{s_1s_2s_1}$to $A_{s_2s_1s_2}$ must map $N$ to $N'$ and $Z$ to
$Z'$, it is not difficult to prove that the algebras
$A_{s_1s_2s_1}$ and $A_{s_2s_1s_2}$ are not isomorphic.

\bigbreak\bigbreak

\end{document}